\documentclass[12pt]{article}
\usepackage{amsmath}
\usepackage{mathrsfs}
\usepackage{amsfonts}
\usepackage{stmaryrd}
\usepackage{latexsym,amssymb,amsmath}
\newcommand{\mbb}{\mathbb}
\newcommand{\ol}{\overline}
\newcommand{\psp}{\vspace{0.2cm}}
\newcommand{\m}{m\!_{_1}}

\begin{document}

\baselineskip=18pt
\newcommand{\headingstobeshown}{}
\def\tenrm{\rm}

\renewcommand{\theequation}{\thesection.\arabic{equation}}


\newenvironment{proof}
               {\begin{sloppypar} \noindent{\it Proof.}}
               {\hspace*{\fill} $\square$ \end{sloppypar}}

\newtheorem{atheorem}{\bf \temp}[section]

\newenvironment{theorem}[1]{\def \temp{#1}
                \begin{atheorem}\medskip
                }
                {\medskip
                 \end{atheorem}}

\newtheorem{thm}[atheorem]{Theorem}
\newtheorem{cor}[atheorem]{Corollary}
\newtheorem{lem}[atheorem]{Lemma}
\newtheorem{pro}[atheorem]{Property}
\newtheorem{prop}[atheorem]{Proposition}
\newtheorem{de}[atheorem]{Definition}
\newtheorem{rem}[atheorem]{Remark}
\newtheorem{fac}[atheorem]{Fact}
\newtheorem{ex}[atheorem]{Example}
\newtheorem{pr}[atheorem]{Problem}
\newtheorem{cla}[atheorem]{Assert}
\newcommand{\la}{\langle}
\newcommand{\ra}{\rangle}
\begin{center}{\Large \bf Structure of Polynomial Representations }\end{center}
\begin{center}{\Large \bf for Orthosymplectic Lie Superalgebras}\end{center}

\vspace{0.2cm}

\begin{center}{\large Cuiling Luo}
\end{center}

\begin{center}{Institute of Mathematics, Academy of Mathematics \&
System Sciences,
\\Chinese Academy of Sciences, Beijing 100190, China}\end{center}
\begin{center}{E-mail: luocuiling@amss.ac.cn}\end{center}
\vspace{0.4cm}
\begin{abstract}

Orthosymplectic Lie superalgebras are fundamental symmetries in
modern physics, such as massive supergravity. However, their
representations are far from being thoroughly understood. In the
present paper, we completely determine the structure of their
various supersymmetric polynomial representations obtained by
swapping bosonic multiplication operators and differential operators
in the canonical supersymmetric polynomial representations. In
particular, we obtain certain new infinite-dimensional irreducible
representations and new composition series of indecomposable
representations for these algebras.

 \vspace{0.5cm} {\bf Keywords:}
representation, Lie superalgebra, weight, composition series.
\end{abstract}
\vspace{0.5cm}
\section{Introduction}

Lie superalgebras were introduced by physicists as the fundamental
tools of studying the supersymmetry in physics (e.g., cf. [2-4],
[7], [26], [29]). For instance, orthosymplectic Lie superalgebras
are symmetries of massive supergravity (cf. [3], [4]). Kac [13] gave
a classification of finite-dimensional Lie superalgebras. Unlike Lie
algebra case, finite-dimensional modules of finite-dimensional
simple Lie superalgebras may not be completely reducible and the
structure of finite-dimensional irreducible modules is much more
complicated due to the existence of so-called {\it atypical} modules
(cf. [14], [15]). Serganova [25] gave a nice survey on  characters
of irreducible representations of simple Lie superalgebras.  Indeed,
representations of orthosymplectic Lie superalgebras are the most
complicated among all the classical Lie superalgebras and people
could so far only get partial information of them.

Palev [24] found the para-Bose and para-Fremi operators as
generators of orthosymplectic Lie superalgebras. Farmer and  Javis
[9] constructed irreducible representations of $osp(3|2)$ and
$osp(4|2)$ by superfield techniques. Moreover, they [10] enumerated
finite-dimensional graded tensor representations of orthosymplectic
Lie superalgebras via standard Young diagrams.  Van der Jeugt [27]
investigated  representations of $osp(3|2)$ by means of the shift
operator technique. Gould and Zhang [12] determined all the
finite-dimensional unitary representations of $osp(2|2n)$. Nishiyama
[21] studied unitary representations of orthosymplectic Lie
superalgebras via supersymettric Heisenberg algebras. He [22] also
obtained the characters and super-characters of discrete series
representations for orthosymplectic Lie superalgebras. Furthermore,
he [23] investigated representations of the superalgebras via super
dual pairs.

Lee Shader [16] investigated certain typical representations of
orthosymplectic Lie superalgebras. Moreover, Benkart, Lee Shader and
Ram [1] studied the tensor product representations of
orthosymplectic Lie superalgebras over the canonical represenations.
Lee Shader [17] obtained certain characteristics of representations
for Lie superalgebras of type C.  Cheng and Zhang [6] found a
combinatorial character formula for orthosymplectic Lie
superalgebras via Howe duality. Dobrev and Zhang [8] classified the
positive energy unitary irreducible representations of superalgebras
$osp(1,2n,\mathbb{R})$. Furthermore,  Cheng, Wang and  Zhang [5]
presented a Fock space approach to representation theory of
$osp(2|2n)$.

Lievens, Stoilova, and Van der Jeugt [18] got  unitary irreducible
representations of the Lie superalgebra $osp(1|2n)$ via the
paraboson Fock space. Moreover, they [19] found a  class of unitary
irreducible representations of the Lie superalegbra $osp(1|2n)$.
Zhang [29] investigated the Schr$\ddot{\mbox{o}}$dinger equation on
the superspace $\mathbb{R}^{m|2n}$ which involved a potential that
varied as an inverse power of the $osp(m|2n)$-invariant distance
from the origin and lead to interesting results regarding the
infinite-dimensional representations of the orthosymplectic Lie
superalgebra $osp(2,m+1|2n)$. In the appendix, he also presented the
structure of a canonical supersymmetric polynomial representation
for $osp(m|2n)$ when $m-2n>1$.

In this paper, we completely determine the structure of  various
supersymmetric polynomial representations of  orthosymplectic Lie
superalgebras obtained by swapping bosonic multiplication operators
and differential operators in the canonical supersymmetric
polynomial representations. In particular, certain new
infinite-dimensional irreducible representations and new composition
series of indecomposable representations for these algebras are
obtained. Below we give a technical introduction to our results.

Denote by $\mbb{Z}$ the ring of integers and by $\mbb{N}$ the set of
nonnegative integers. For convenience, we also use the following
notation of indices:
\begin{equation}
\ol{i,j}=\{i,i+1,...,j\},\end{equation} where $i\leq j$ are
integers. Moreover, we also use $\{0,1\}$ to denote
$\mbb{Z}_2=\mbb{Z}/2\mbb{Z}$ when the context is clear. Let
$E_{i,j}$ be the square matrix whose $(i,j)$-entry is 1 and the
others are zero. The general linear Lie superalgebra of
$(m+2n)\times (m+2n)$ matrices $gl(m,2n)=gl(m,2n)_0\oplus
gl(m,2n)_1$ with
\begin{equation}
gl(m,2n)_0=\sum_{i,j=1}^{m}\mbb{C}E_{i,j}+\sum_{s,t=1}^{2n}\mbb{C}E_{m+s,m+t},\;
gl(m,2n)_1=\sum_{i=1}^{m}\sum_{s=1}^{2n}(\mbb{C}E_{i,m+s}+\mbb{C}E_{m+s,i})\end{equation}
and the Lie superbracket:
\begin{equation}
[u,v]=uv-(-1)^{\iota_{_1}\iota_{_2}}vu\qquad\mbox{for}\;u\in
gl(m,2n)_{\iota_{_1}},\;v\in gl(m,2n)_{\iota_{_2}}.\end{equation}
Assume that $m=2m_1$ is an even integer. The orthosymplectic Lie
superalgebra $osp(m,2n)$ is the subalgebra of $gl(m,2n)$ consisting
of the matrices of the form
\begin{equation}
\left(
\begin{array}{cccc}
A & B & H & H_1\\
C &-A^T & K & K_1\\
K_1^T & H_1^T &D & E\\
-K^T &-H^T & F & -D^T
\end{array}
\right)
\end{equation}
where $A$, $B$ and $C$ are $m_1\times m_1$ matrices such that
$B=-B^T$, $C=-C^T$; $D$, $E$ and $F$ are $n\times n$ matrices such
that $E=E^T$, $F=F^T$; $H$, $H_1$, $K$ and $K_1$ are $m_1\times n$
matrices. Let ${\cal
A}=\mathbb{C}[x_1,\cdots,x_m,\theta_1,\cdots,\theta_{2n}]$ be the
polynomial algebra in bosonic variables $x_1,\cdots,x_m$ and
fermionic variables $\theta_1,\cdots,\theta_{2n}$, i.e.
\begin{equation}
x_ix_j=x_jx_i,\;\;\theta_p\theta_q=-\theta_q\theta_p,\;\;x_i\theta_p=\theta_px_i,\;\;i,j\in\overline{1,m},\;\;p,q\in\overline{1,2n}.
\end{equation}
 Taking $r\in \ol{0,m}$, we have the following
 supersymmetric polynomial representation of $gl(m|2n)$:
\begin{equation}
\begin{array}{l}
E_{i,j}|_{\mathcal{A}}=\left\{
\begin{array}{lll}
-x_j\partial_{x_i}-\delta_{i,j} &\mbox{if}&i,j\in\overline{1,r},\\
\partial_{x_i}\partial_{x_j} &\mbox{if} &i\in\overline{1,r},\;j\in\overline{r+1,m},\\
-x_ix_j &\mbox{if} &j\in\overline{1,r},\;i\in\overline{r+1,m},\\
x_i\partial_{x_j}& \mbox{if} &i,j\in\overline{r+1,m},
\end{array}\right.\\
\\
E_{i,m+p}|_{\mathcal{A}}=\left\{
\begin{array}{lll}
\partial_{x_i}\partial_{\theta_p} &\mbox{if} &i\in\overline{1,r},\;p\in\overline{1,2n},\\
x_i\partial_{\theta_p} &\mbox{if}
&i\in\overline{r+1,m},\;p\in\overline{1,2n},
\end{array}\right.\\
\\
E_{m+p,j}|_{\mathcal{A}}=\left\{
\begin{array}{lll}
-\theta_px_j &\mbox{if} &j\in\overline{1,r},\;p\in\overline{1,2n}\\
\theta_p\partial_{x_j} &\mbox{if}&
j\in\overline{r+1,m},\;p\in\overline{1,2n},
\end{array}
\right.\\
\\
E_{m+p,m+q}|_{\mathcal{A}}=\theta_p\partial_{\theta_q}\;\mbox{if
}p,q\in\overline{1,2n}, \label{I1}
\end{array}
\end{equation}
which is obtained from the canonical supersymmetric polynomial
representation ($r=0$) by swapping $\partial_{x_i}$ and $-x_i$ for
$i\in\overline{1,r}$. In particular, we have the restricted
representation of $osp(m,2n)$ on ${\cal A}$.

For $k\in\mbb{N}$, we define
\begin{eqnarray}
&&{\cal
A}_k^{r}=\mbox{Span}\;\{x^\alpha\theta_{i_1}\cdots\theta_{i_t}\in{\cal
A}^r\mid 0\leq t\leq
2n;\;i_1,\cdots,i_t\in\overline{1,2n};\;\alpha\in\mathbb{N}^{m};\nonumber\\
&&t-\sum\limits_{i=1}^r\alpha_i+\sum\limits_{j=r+1}^{m}\alpha_j=k\}.\label{I3}
\end{eqnarray}
Then subspaces ${\cal A}^{r}_k$ form $osp(m,2n)$-submodules. We can
assume $r\leq m_1$ by symmestry. Denote
\begin{equation}\Delta=-\sum\limits_{i=1}^r
x_i\partial_{x_{m_1+i}}+\sum\limits_{i=r+1}^{m_1}\partial_{x_i}\partial_{x_{m_1+i}}+
\sum\limits_{j=1}^n\partial_{\theta_j}\partial_{\theta_{n+j}},\end{equation}
\begin{equation}\eta=
\sum\limits_{i=1}^rx_{m_1+i}\partial_{x_i}+\sum\limits_{i=r+1}^{\m}x_ix_{m_1+i}+\sum\limits_{j=1}^n\theta_j\theta_{n+j}.
\end{equation}
Set
\begin{equation}
{\cal H}^r_k=\{f\in{\cal A}^r_k\mid\Delta(f)=0\}.
\end{equation}
Denote by $\langle F\rangle$ the $osp(m,2n)$-submodule generated by
a subset $F$.\psp

\textbf{Theorem 1.} {\it We have the following conclusions:

1) Assume $r=0$. If
 $k>2(n-m_1+1)$ or $k\leq (n-m_1+1)$, then ${\cal A}^{0}_k={\cal
H}^{0}_k\oplus\eta{\cal A}^{0}_{k-2}$ and the subspace ${\cal
H}^{0}_k$ is an irreducible $osp(m,2n)$-submodule. When
$(n-m_1+1)<k\leq 2(n-m_1+1)$,
\begin{equation}
{\cal H}_k^{0}\supset\eta^{k-(n-m_1+1)}{\cal
H}^{0}_{2(n-m_1+1)-k}\supset\{0\}
\end{equation}
is a composition series.

2) Suppose $r=m_1$. We always have ${\cal A}^{\m}_k={\cal
H}^{\m}_k\oplus\eta{\cal A}^{\m}_{k-2}$. If $k>n$, ${\cal
H}^{\m}_k=\{0\}$. When $k\leq n$, the subspace ${\cal H}^{m_1}_k$ is
an irreducible $osp(m,2n)$-submodule.

3) Assume $0<r<m_1$. If $k\leq n-m_1+r+1$, then ${\cal A}^r_k={\cal
H}^r_k\bigoplus\eta{\cal A}^r_{k-2}$ and the submodule ${\cal
H}^r_k$ is irreducible. When $k>n-m_1+r+1$, we have the following
composition series
\begin{eqnarray}
&&{\cal H}^r_k\supset\eta^{k-n+m_1-r-1}{\cal
H}^r_{-k+2(n-m_1+r+1)}\supset\{0\}\mbox{ if }r<m_1-1;\label{d6252}\\
&&{\cal H}^{m_1-1}_k\supset\la x_{m_1}^k\ra\supset\eta^{k-n}{\cal
H}^{m_1-1}_{-k+2n}\supset\{0\}.\end{eqnarray}} \psp

Take a subset $T$ of $\overline{1,2n}$. Denote
$\bar{T}=\overline{1,2n}\setminus T$. Let $osp(m,2n)$ act on ${\cal
A}'=\mathbb{C}[x_1,\cdots,x_{2n};\theta_1,\cdots,\theta_{m}]$ via
\begin{equation}
\begin{array}{l}
E_{p,q}|_{\mathcal{A}'}=\theta_p\partial_{\theta_q}\qquad\mbox{ if }p,q\in\overline{1,m},\\
\\
E_{m+i,m+j}|_{\mathcal{A}'}=\left\{
\begin{array}{lll}
-x_j\partial_{x_i}-\delta_{i,j} &\mbox{if}&i,j\in T,\\
\partial_{x_i}\partial_{x_j} &\mbox{if} &i\in T,\;j\in \bar{T},\\
-x_ix_j &\mbox{if} &i\in\bar{T},\;j\in T,\\
x_i\partial_{x_j}& \mbox{if} &i,j\in\bar{T},
\end{array}\right.\\
\\
E_{m+i,p}|_{\mathcal{A}'}=\left\{
\begin{array}{lll}
\partial_{x_i}\partial_{\theta_p} &\mbox{if} &i\in T,\;p\in\overline{1,m},\\
x_i\partial_{\theta_p} &\mbox{if} &i\in\bar{T},\;p\in\overline{1,m},
\end{array}\right.\\
\\
E_{p,m+j}|_{\mathcal{A}'}=\left\{
\begin{array}{lll}
-\theta_px_j &\mbox{if} &j\in T,\;p\in\overline{1,m},\\
\theta_p\partial_{x_j} &\mbox{if}& j\in\bar{T},\;p\in\overline{1,m}.
\end{array}
\right.
\end{array}\label{I6181}
\end{equation}
Then we obtain another supersymmetric polynomial representation of
$osp(m,2n)$, which is obtained from the corresponding canonical one
($T=\emptyset$) by swapping $\partial_{x_i}$ and $-x_i$ for $i\in
T$. The subspace
\begin{eqnarray}
&&{\cal
A}'_k=\mbox{Span}\:\{x^\alpha\theta_{i_1}\cdots\theta_{i_t}\in{\cal
A}'\mid t\in\overline{0,m};\;
i_1,\cdots,i_t\in\overline{1,m};\;\alpha\in\mathbb{N}^{2n};\nonumber\\
&&\sum\limits_{i\in\bar{T}}\alpha_i-\sum\limits_{i\in
T}\alpha_i=k-t\}\label{i61}
\end{eqnarray}
 forms an $osp(m,2n)$-submodule for $k\in\mbb{N}$. Denote
\begin{equation}
S_1=\{i\in\overline{1,2n}\mid
i\in\bar{T},\;n+i\in\bar{T}\},\;\;T_1=\{i\in\overline{1,2n}\mid i\in
T,\;n+i\in T\}.
\end{equation}

\textbf{Theorem 2.} {\it The following statements hold:

1) The submodule ${\cal A}'_k$ is irreducible when $S_1\cup
T_1\neq\emptyset$. In particular, ${\cal A}'_k$ is not highest
weight type if $S_1\neq\emptyset$ and $T_1\neq\emptyset$.

2) Suppose $S_1=\emptyset$ and $T_1=\emptyset$. We may assume
$T=\overline{1,n}$ by symmetry.

 a) The submodule ${\cal A}'_k$ is
irreducible when $k\neq m_1$.

 b) The submodule ${\cal A}'_{\m}=\la
\theta_1\cdots\theta_{\m}\ra\oplus\la
(x_{n-1}x_{2n}-x_nx_{2n-1})\theta_1\cdots\theta_{\m}\ra$ is a direct
sum of two irreducible submodules.} \psp

Suppose that $m=2m_1+1$ is an odd integer. The orthosymplectic Lie
superalgebra $osp(m,2n)$ is the subalgebra of $gl(m,2n)$ consisting
of the matrices of the form
\begin{equation}
\left(
\begin{array}{ccccc}
A & B & U & H & H_1\\
C &-A^T &V& K & K_1\\
-V^T& U^T&0&M& M_1\\
K_1^T & H_1^T& M_1^T &D & E\\
-K^T &-H^T& -M^T & F & -D^T
\end{array}
\right)
\end{equation}
where $A$, $B$ and $C$ are $m_1\times m_1$ matrices such that
$B=-B^T$, $C=-C^T$; $D$, $E$ and $F$ are $n\times n$ matrices such
that $E=E^T$, $F=F^T$; $H$, $H_1$, $K$ and $K_1$ are $m_1\times n$
matrices; $U$ and $V$ are $m_1\times1$ matrices; $M$ and $M_1$ are
$1\times n$ matrices. Similarly, we have a representation of
$osp(m,2n)$ on ${\cal
A}=\mathbb{C}[x_1,\cdots,x_m,\theta_1,\cdots,\theta_{2n}]$ via
(1.6), (1.17) and a representation of $osp(m,2n)$ on ${\cal
A}'=\mathbb{C}[x_1,\cdots,x_{2n};\theta_1,\cdots,\theta_{m}]$ via
(1.14), (1.17). \psp

\textbf{Theorem 3.} {\it All the $osp(2m+1,2n)$-submodules ${\cal
H}^r_k$ and ${\cal A}'_k$ are irreducible.} \psp

In addition to the results given in the above, we have also
constructed a basis for the module ${\cal H}^r_k$, and the
submodules $\la \theta_1\cdots\theta_{\m}\ra$ and $\la
(x_{n-1}x_{2n}-x_nx_{2n-1})\theta_1\cdots\theta_{\m}\ra$ in b) of 2)
in Theorem 2.

In Section 2, we prove Theorem 1. Moreover, Theorem 2 is proved in
Section 3. We give a proof of Theorem 3 in Section 4.

\section{Proof of Theorem 1}

\setcounter{equation}{0} In this section, we discuss the polynomial
representations of $osp(2m_1,2n)$ ($m_1>0,n>0$) defined via (1.6).

Recall the Lie superalgebra $osp(2m_1,2n)$ given (1.4). The even
part of $osp(2m,2n)$
\begin{eqnarray}
&&osp(2m_1,2n)_{0}=\sum\limits_{i,j=1}^{\m}\big[\mathbb{C}(E_{i,j}-E_{\m+j,\m+i})+\mathbb{C}(E_{i,\m+j}-E_{j,\m+i})\nonumber\\
&&+\mathbb{C}(E_{\m+i,j}-E_{\m+j,i})\big]
+\sum\limits_{p,q=1}^n\big[\mathbb{C}(E_{2\m+p,2\m+q}-E_{2\m+n+q,2\m+n+p})\nonumber\\
&&+\mathbb{C}(E_{2\m+p,2\m+n+q}+E_{2\m+q,2\m+n+p})\nonumber\\
& &+\mathbb{C}(E_{2\m+n+p,2\m+q}+E_{2\m+n+q,2\m+p})\big]
\end{eqnarray}
is a subalgebra isomorphic to $o(2m_1,\mbb{C})\oplus sp(2n,\mbb{C})$
and the odd part is
\begin{eqnarray}
osp(2m_1,2n)_1&=&\sum_{i=1}^{\m}\sum_{p=1}^n\big[\mathbb{C}(E_{i,2\m+p}-E_{2\m+n+p,\m+i})+\mathbb{C}(E_{i,2\m+n+p}+E_{2\m+p,\m+i})\nonumber\\
&&+\mathbb{C}(E_{\m+i,2\m+p}-E_{2\m+n+p,i})+\mathbb{C}(E_{\m+i,2\m+n+p}+E_{2\m+p,i})\big].
\end{eqnarray}
Take
\begin{equation}
H=\sum\limits_{i=1}^{\m}\mathbb{C}(E_{i,i}-E_{\m+i,\m+i})+\sum\limits_{j=1}^n\mathbb{C}(E_{2\m+j,2\m+j}-E_{2\m+n+j,2\m+n+j})
\end{equation}
as a Cartan subalgebra of $osp(2\m,2n)$. Let
$\lambda_1,\cdots,\lambda_{\m}$, $\nu_1,\cdots,\nu_n$ be the
fundamental weights of $o(2m_1,\mbb{C})\oplus sp(2n,\mbb{C})$. Let
\begin{eqnarray}
osp(2\m,2n)^+&=&\sum\limits_{1\leq i<j\leq
\m}\big(\mathbb{C}(E_{i,j}-E_{\m+j,\m+i})+\mathbb{C}(E_{i,\m+j}-E_{j,\m+i})\big)\nonumber\\
&&+\sum\limits_{1\leq p<q\leq
n}\mathbb{C}(E_{2\m+p,2\m+q}-E_{2\m+n+q,2\m+n+p})\nonumber
\\&&+\sum\limits_{1\leq p\leq q\leq
n}\mathbb{C}(E_{2\m+p,2\m+n+q}+E_{2\m+q,2\m+n+p})\nonumber\\
&&+\sum\limits_{1\leq i\leq \m,1\leq q\leq
n}\big(\mathbb{C}(E_{i,2\m+q}-E_{2\m+n+q,\m+i})\nonumber\\
&&+\mathbb{C}(E_{i,2\m+n+q}+E_{2\m+q,\m+i})\big),
\end{eqnarray}
and $osp(2\m,2n)^+_{\sigma}=osp(2\m,2n)_{\sigma}\cap osp(2\m,2n)^+$
for $\sigma=0,1$. A weight vector $f\in{\cal A}^r$ is called a
\textit{highest weight vector} if $osp(2\m,2n)^+(f)=0$, and the
corresponding weight is called the \textit{highest weight}.

We first quote a useful lemma found by Xu [\ref{r11}].
\begin{lem}
Suppose ${\cal A}$ is a free module over a subalgebra $B$ generated
by a filtrated subspace $V=\bigcup_{r=0}^\infty V_r$ (i.e.,
$V_r\subset V_{r+1}$). Let ${\cal T}_1$ be a linear operator on
${\cal A}$ with a right inverse ${\cal T}_1^-$ such that
\begin{equation}
{\cal T}_1(B),{\cal T}_1^-(B)\subset B,\ \ {\cal
T}_1(\eta_1\eta_2)={\cal T}_1(\eta_1)\eta_2,\ \ {\cal
T}_1^-(\eta_1\eta_2)={\cal T}_1^-(\eta_1)\eta_2
\end{equation}
for $\eta_1\in B$, $\eta_2\in V$, and let ${\cal T}_2$ be a linear
operator on ${\cal A}$ such that
\begin{equation}
{\cal T}_2(V_{r+1})\subset BV_r,\ \ {\cal T}_2(f\zeta)=f{\cal
T}_2(\zeta)\ \ for\ \ 0\leq r\in\mathbb{Z},\ \ f\in B,\ \ \zeta\in
{\cal A}.
\end{equation}
Then we have
\begin{equation}
\begin{array}{ll}
&\{g\in {\cal A}\mid ({\cal T}_1+{\cal T}_2)(g)=0\}\\
=&\mbox{\it Span}\:\{\sum\limits_{i=0}^\infty(-{\cal T}_1^-{\cal
T}_2)^i(hg)|g\in V,h\in B;{\cal T}_1(h)=0\},
\end{array}
\end{equation}
where the summation is finite under our assumption.
\end{lem}

Below we discuss the $osp(2\m,2n)$-module structure of ${\cal
A}^r_k$ case by case. \psp

\textbf{Case 1}. $r=0$

In this case,
\begin{equation}
\Delta=\sum\limits_{i=1}^{\m}\partial_{x_i}\partial_{x_{\m+i}}+\sum\limits_{j=1}^n\partial_{\theta_j}\partial_{\theta_{n+j}},\;\;
\eta=\sum\limits_{i=1}^{\m}x_ix_{\m+i}+\sum\limits_{j=1}^n\theta_j\theta_{n+j}.\end{equation}
The subspaces ${\cal A}^0_k$ ($k\in\mathbb{N}$) are all finite
dimensional and ${\cal A}^0_k=0$ when $k<0$. Moreover, ${\cal
A}^0_k={\cal H}^0_k\oplus\eta{\cal A}_{k-2}^0$ when $\m=0$ or $n=0$.
\begin{thm}
If $k>2(n-\m+1)$ or $k\leq (n-\m+1)$, then ${\cal A}^{0}_k={\cal
H}^{0}_k\oplus\eta{\cal A}^{0}_{k-2}$ and the subspace ${\cal
H}^{0}_k$ is an irreducible $osp(2\m,2n)$-submodule with highest
weight vector $x_1^k$ and the corresponding highest weight
$k\lambda_1$. When $(n-\m+1)<k\leq 2(n-\m+1)$,
\begin{equation}
{\cal H}_k^{0}\supset\eta^{k-(n-\m+1)}{\cal
H}^{0}_{2(n-\m+1)-k}\supset\{0\}
\end{equation}
is a composition series. Moreover,
\begin{eqnarray}
&&\{\sum\limits_{\stackrel{r_1,\cdots,r_{\m-1}\in\mathbb{N},}{
l_1,\cdots,l_n\in\{0,1\}}}\frac{(-1)^{\sum\limits_{i=1}^{\m-1}r_i+\sum\limits_{j=1}^nl_j}\alpha_{\m}!\alpha_{2\m}!
\prod\limits_{i=1}^{\m-1}r_i!\left(\begin{array}{c}\alpha_i\\
r_i\end{array}\right)\left(\begin{array}{c}\alpha_{\m+i}\\
r_i\end{array}\right)\prod\limits_{j=1}^n\delta_{l_j,1}(-1)^{\beta_j}\beta_j\beta_{n+j}}{(\alpha_{\m}+\sum\limits_{i=1}^{\m-1}r_i+\sum\limits_{j=1}^nl_j)!(\alpha_{2\m}+\sum\limits_{i=1}^{\m-1}r_i+\sum\limits_{j=1}^nl_j)!}\nonumber\\
&&\times
x_{\m}^{\alpha_{\m}+\sum\limits_{i=1}^{\m-1}r_i+\sum\limits_{j=1}^nl_j}
x_{2\m}^{\alpha_{2\m}+\sum\limits_{i=1}^{\m-1}r_i+\sum\limits_{j=1}^nl_j}\prod\limits_{i=1}^{\m-1}x_i^{\alpha_i-r_i}x_{\m+i}^{\alpha_{\m+i}-r_i}\prod\limits_{j=1}^n\theta_j^{\beta_j-l_j}\beta_{n+j}^{\beta_{n+j}-l_j}\nonumber\\
&&\mid \alpha_1,\cdots,\alpha_{2\m}\in\mathbb{N},
l_1\cdots,l_{2n}\in\{0,1\};\;\alpha_{\m}\alpha_{2\m}=0;\;\sum\limits_{i=1}^{2\m}\alpha_i+\sum\limits_{j=1}^{2n}l_j=k
\}
\end{eqnarray}
is a basis of ${\cal H}^{0}_k\;(k>0)$.
\end{thm}
\begin{proof} We divide our arguments as the following steps.\psp

 (1) \textsl{ The submodule ${\cal H}^{0}_k$ is generated by
$x_1^k$ for $k>0$.}

\psp It is well known that ${\cal H}^0_k=\la x_1^k\ra$ when $n=0$.
Now assume $k\geq2$. Take induction on $n$. For any $0\neq f\in{\cal
H}^{0}_k$, we can write
\begin{equation}
f=f_0+f_1\theta_n+f_2\theta_{2n}+f_3\theta_n\theta_{2n},
\end{equation}
with
\begin{equation}f_i\in\mathbb{C}[x_1,\cdots,x_{2\m};\theta_1\cdots,\theta_{n-1},\theta_{n+1},\cdots,\theta_{2n-1}].
\end{equation}
We denote
\begin{equation}\Delta'=\sum\limits_{i=1}^{\m}\partial_{x_i}\partial_{x_{\m+i}}
+\sum\limits_{j=1}^{n-1}\partial_{\theta_j}\partial_{\theta_{n+j}}\end{equation}
and get
\begin{equation}
0=\Delta(f)=\Delta'f_0+(\Delta'f_1)\theta_n+(\Delta'f_2)\theta_{2n}+(\Delta'f_3)\theta_n\theta_{2n}-f_3.
\end{equation}
Hence
\begin{equation}
\Delta'(f_1)=\Delta'(f_2)=\Delta'(f_3)=0,\;\;\Delta'(f_0)-f_3=0.
\end{equation}
By induction,
\begin{equation}f_1=X'(x_1^{k-1}),\;\;f_2=X''(x_1^{k-1}),\;\;
f_3=X(x_1^{k-2})\end{equation} for some $X',X'',X\in
U(osp(2\m,2(n-1)))$. Thus
\begin{equation}
f_1\theta_n=\frac{1}{k}X'(E_{\m+1,2\m+2n}+E_{2\m+n,1})(x_1^k),\end{equation}
\begin{equation}f_2\theta_{2n}=\frac{1}{k}X''(E_{2\m+2n,1}-E_{\m+1,2\m+n})(x_1^k).
\end{equation}
Moreover,
\begin{equation}f_3=X(x_1^{k-2})=\Delta'\big(\frac{1}{k-1}X(x_1^{k-1}x_{\m+1})\big)\end{equation}
because $[\Delta',X]=0$. So
\begin{eqnarray}
0&=&\Delta'(f_0)-f_3=\Delta'(f_0)-\Delta'\big(\frac{1}{k-1}X(x_1^{k-1}x_{\m+1})\big)\nonumber\\
&=& \Delta'(f_0-\frac{1}{k-1}X(x_1^{k-1}x_{\m+1})),
\end{eqnarray}
which implies
\begin{equation}f_0-\frac{1}{k-1}X(x_1^{k-1}x_{\m+1})\in
U(osp(2\m,2(n-1)))(x_1^k).\end{equation} Note
\begin{eqnarray}
&&\frac{1}{k-1}x_1^{k-1}x_{\m+1}+x_1^{k-2}\theta_n\theta_{2n}\nonumber\\
=&&\frac{1}{k(k-1)}(E_{\m+1,2\m+n}-E_{2\m+2n,m})(E_{\m+1,2\m+2n}+E_{2\m+n,1})(x_1^k).
\end{eqnarray}
Thus
\begin{eqnarray}
f&=&(f_0-\frac{1}{k-1}X(x_1^{k-1}x_{\m+1}))+f_1\theta_n+f_2\theta_{2n}
\nonumber\\
&&+X(\frac{1}{k-1}x_1^{k-1}x_{\m+1}+x_1^{k-2}\theta_n\theta_{2n})\in\la
x_1^k\ra.\end{eqnarray}

(2) \textsl{The submodule ${\cal H}^{0}_k$ is irreducible if
$k\leq(n-\m+1)$ or $k>2(n-\m+1)$.}

\psp

Denote
\begin{equation}\eta_x=
\sum\limits_{i=1}^{\m}x_ix_{\m+i},\;\;
\eta_\theta=\sum\limits_{j=1}^n\theta_j\theta_{n+j}.
\end{equation}
Then $\eta=\eta_x+\eta_\theta$. By direct calculation, we get all
the weight vectors in ${\cal A}^0_k$ annihilated by $osp(2m,2n)^+_0$
are  scalar multiples of the following elements
\begin{equation}\sum\limits_{i=0}^l
a_ix_1^{k-2l-t}\eta_x^{l-i}\eta_\theta^i\theta_1\cdots\theta_t
\end{equation}
 with
$l-2l-t\geq0$, $l\geq0$, $0\leq t\leq n$ and $a_i=0$ for $i>n-t$.
Suppose $\sum\limits_{i=0}^l
a_ix_1^{k-2l-t}\eta_x^{l-i}\eta_\theta^i\theta_1\cdots\theta_t\in{\cal
H}^0_k$. Then
\begin{eqnarray}
0&=&\Delta(\sum\limits_{i=0}^l
a_ix_1^{k-2l-t}\eta_x^{l-i}\eta_\theta^i\theta_1\cdots\theta_t)\nonumber\\
&=&\sum\limits_{i=0}^{l-1}
a_i(l-i)(k-l-t+\m-1-i)x_1^{k-2l-t}\eta_x^{l-i-1}\eta_\theta^i\theta_1\cdots\theta_t\nonumber\\
&&-\sum\limits_{i=1}^l
a_ii(n-t-i+1)x_1^{k-2l-t}\eta_x^{l-i}\eta_\theta^{i-1}\theta_1\cdots\theta_t
\end{eqnarray}
by (4.31) in [28]. So
\begin{equation}a_i(l-i)(k-l-t+\m-1-i)-a_{i+1}(i+1)(n-t-i)=0\end{equation} for $0\leq i<n-t$
and $l\leq n-t$. Thus up to a scalar multiple,  all the weight
vectors in ${\cal H}^0_k$ annihilated by $osp(2\m,2n)^+_{0}$ are
\begin{equation}
f_{l,t}=\sum\limits_{i=0}^l
\frac{(n-t-i)!}{i!(l-i)!(k-l-t+\m-1-i)!}x_1^{k-2l-t}\eta_x^{l-i}\eta_\theta^i\theta_1\cdots\theta_t,\label{d5}
\end{equation}
for $0\leq t\leq n$ and $0\leq
l\leq\mbox{min}\{n-t,\frac{1}{2}(k-t)\}$.

Note
\begin{equation}
f_{0,t-1}=(-1)^{t-1}\frac{n-t+1}{k-t+\m}(E_{1,2\m+t}-E_{2\m+n+t,\m+1})(f_{0,t})\label{d1103}
\end{equation}
for $0<t\leq n$ and
\begin{equation}
(E_{1,2\m+n+t+1}+E_{t+1,\m+1})(f_{l,t})=(-1)^{t-1}(k+\m-n-1-l)f_{l-1,t+1}\label{d6}
\end{equation}
for $0<l\leq\mbox{min}\{n-t,\frac{1}{2}(k-t)\}$.  If $k\leq
(n-\m+1)$, then
\begin{equation}
k+\m-n-1-l\leq-l<0.
\end{equation}
When $k>2(n-\m+1)$,
\begin{equation}
k+\m-n-1-l\geq k+\m-n-1-\frac{1}{2}(k-t)\geq
\frac{1}{2}k-(n-\m+1)>0.
\end{equation}
Thus $x_1^k$ is the only highest weight vector in ${\cal H}^0_k$ up
to a scalar multiple, which implies ${\cal H}^0_k$ is irreducible by
(1). \psp

(3) \textsl{${\cal A}_k^{0}={\cal H}^{0}_k\oplus\eta{\cal
A}^{0}_{k-2}$ when $k\leq(n-\m+1)$ or $k>2(n-\m+1)$.}

\psp Since $x_1^k\notin\eta{\cal A}^{0}_{k-2}$, we get ${\cal
H}^{0}_k\cap\eta{\cal A}^{0}_{k-2}=0$. Now we still have to show
that ${\cal A}^0_{k}={\cal H}^0_{k}+\eta{\cal A}^0_{k-2}$. Note that
${\cal A}^0_{k}$ is a finite-dimensional $osp(2\m,2n)_0$-module in
this case. It is sufficient to check
\begin{equation}
x_1^{k-2l-t}\theta_1\cdots\theta_t\eta_x^{l-p}\eta_\theta^p\in{\cal
H}^{0}_k+\eta{\cal A}^{0}_{k-2}\;\;\mbox{for all}\; 0\leq t\leq
n;\;p\leq l\in\mathbb{N}.\label{d3}
\end{equation}
Observe
\begin{equation}
x_1^{k-2l-t}\theta_1\cdots\theta_t\eta_x^{l-p}\eta_\theta^p=\eta
\sum\limits_{i=0}^{n-t-p}(-1)^ix_1^{k-2l-t}\theta_1\cdots\theta_t\eta_x^{l-p-i-1}\eta_\theta^{p+i}
\end{equation}
if $l>n-t$. When $l\leq n-t$, we have
\begin{equation}
\Delta\big(\sum\limits_{s=0}^l\frac{(n-t-s)!}{s!(l-s)!(k-l-t+n-s-1)!}
x_1^{k-2l-t}\theta_1\cdots\theta_t\eta_x^{l-s}\eta_\theta^s\big)=0,
\end{equation}
\begin{eqnarray}
x_1^{k-2l-t}\theta_1\cdots\theta_t\eta_x^{l-p}\eta_\theta^p&=&\eta
\big(\sum\limits_{i=0}(-1)^ix_1^{k-2l-t}\theta_1\cdots\theta_t\eta_x^{l-p-i-1}\eta_\theta^{p+i}\big)\nonumber\\
&&+(-1)^{s-p}x_1^{k-2l-t}\theta_1\cdots\theta_t\eta_x^{l-s}\eta_\theta^s
\end{eqnarray}
for $s>p$ and
\begin{eqnarray}
x_1^{k-2l-t}\theta_1\cdots\theta_t\eta_x^{l-p}\eta_\theta^p&=&\eta
\big(\sum\limits_{i=0}^{p-s-1}(-1)^ix_1^{k-2l-t}\theta_1\cdots\theta_t\eta_x^{l-p+i}\eta_\theta^{p-i-1}\big)\nonumber\\
&&+(-1)^{p-s}x_1^{k-2l-t}\theta_1\cdots\theta_t\eta_x^{l-s}\eta_\theta^s
\end{eqnarray}
for $s<p$. Thus
\begin{eqnarray}
&&\big(\sum\limits_{s=0}^l\frac{(-1)^{p-s}(n-t-s)!}{s!(l-s)!(k-l-t+n-s-1)!}\big)
x_1^{k-2l-t}\theta_1\cdots\theta_t\eta_x^{l-p}\eta_\theta^p\nonumber\\
=&&\sum\limits_{s=0}^l\frac{(n-t-s)!}{s!(l-s)!(k-l-t+n-s-1)!}
x_1^{k-2l-t}\theta_1\cdots\theta_t\eta_x^{l-s}\eta_\theta^s\nonumber\\
&&+\eta\big(\sum\limits_{s=0}^{p-1}\sum\limits_{i=0}^{p-s-1}(-1)^{p-s+i}x_1^{k-2l-t}\theta_1\cdots\theta_t\eta_x^{l-p+i}\eta_\theta^{p-i-1}\nonumber\\
&&+\sum\limits_{s=p+1}^l\sum\limits_{i=0}(-1)^{i+s-p}x_1^{k-2l-t}\theta_1\cdots\theta_t\eta_x^{l-p-i-1}\eta_\theta^{p+i}\big),
\end{eqnarray}
which implies (\ref{d3}) holds.\psp

(4) \textsl{${\cal H}^{0}_k$ has only one nonzero proper submodule
$\eta^{k-(n-\m+1)}{\cal H}^{0}_{2(n,\m+1)-k}$ when $(n-\m+1)<k\leq
2(n-\m+1)$.}\psp

According to (2.28), we have
\begin{equation}
{\cal H}_k^{0}=\bigoplus\limits_{\stackrel{0\leq t\leq n,}{0\leq
l\leq\rm{min}\{n-t,\frac{1}{2}(k-t)\}}}
U(osp(2\m,2n)_{0}(f_{l,t}).\end{equation} Since
\begin{equation}
f_{k-(n-\m+1),0}=\frac{1}{(k-n+\m-1)!}\eta^{k-n+\m-1}x_1^{2(n-\m+1)-k}\in{\cal
H}^0_k
\end{equation}
and ${\cal H}^0_{2(n-\m+1)-k}=\la x_1^{2(n-\m+1)-k}\ra$, we have
\begin{equation}
\eta^{k-(n-\m+1)}{\cal H}^0_{2(n-\m+1)-k}=\eta^{k-(n-\m+1)}\la
x_1^{2(n-\m+1)-k}\ra\subset{\cal H}^0_k.
\end{equation}
Note
\begin{equation}
\Delta\big(\sum\limits_{i=0}^{l-k+(n-\m+1)}
\frac{(n-i)!\eta_x^{l-k+n-\m+1-i}\eta_\theta^i}{i!(l-k+n-\m+1-i)!(n-l-i)!}x_1^{k-2l}\big)=0,
\end{equation}
which implies
\begin{eqnarray}
&& \eta^{k-(n-\m+1)} (\sum\limits_{i=0}^{l-k+(n-\m+1)}
\frac{(n-i)!\eta_x^{l-k+n-\m+1-i}\eta_\theta^i}{i!(l-k+n-\m+1-i)!(n-l-i)!}x_1^{k-2l})\nonumber\\
&&\in\eta^{k-(n-\m+1)}\la x_1^{2(n-\m+1)-k}\ra\subset{\cal H}^0_k
\end{eqnarray}
for $l\geq k-(n-\m+1)$. Observe
\begin{equation}
osp(2\m,2n)_{0}^+\eta^{k-(n-\m+1)} (\sum\limits_{i=0}^{l-k+(n-\m+1)}
\frac{(n-i)!\eta_x^{l-k+n-\m+1-i}\eta_\theta^i}{i!(l-k+n-\m+1-i)!(n-l-i)!}x_1^{k-2l})=0.
\end{equation}
Thus we get
\begin{eqnarray}
\eta^{k-(n-\m+1)}\big(\sum\limits_{i=0}^{l-k+(n-\m+1)}
\frac{(n-i)!\eta_x^{l-k+n-\m+1-i}\eta_\theta^i}{i!(l-k+n-\m+1-i)!(n-l-i)!}x_1^{k-2l}\big)=c_lf_{l,0}
\end{eqnarray}
for some $c_l\in\mathbb{C}$ by (\ref{d5}). Hence
$f_{l,0}\in\eta^{k-(n-\m+1)}{\cal H}^{0}_{2(n-\m+1)-k}$.
Consequently, $f_{l,t}\in\eta^{k-(n-\m+1)}{\cal
H}^{0}_{2(n-\m+1)-k}$ by (\ref{d6}) for all $l>k-(n-\m+1)$.

Suppose that $W$ is a nonzero submodule of ${\cal H}^{0}_k$. If
there exists $f_{l,t}\in W$ for some $l<k-(n-\m+1)$, then
$f_{0,l+t}\in W$ by (\ref{d6}), which implies $x_1^k\in W$ by
(\ref{d1103}); otherwise $W\subset\eta^{k-(n-\m+1)}{\cal
H}^{0}_{2(n-\m+1)-k}$, which means $W=\eta^{k-(n-\m+1)}{\cal
H}^{0}_{2(n-\m+1)-k}$ since $\eta^{k-(n-\m+1)}{\cal
H}^{0}_{2(n-\m+1)-k}$ is an irreducible $osp(2\m,2n)$-submodule.
\end{proof}
\psp

{\bf Remark}. In the appendix of [29], Zhang presented the structure
of a canonical supersymmetric polynomial representation for
$osp(m|2n)$ when $m-2n>1$. Our above theorem give a complete answer
to the structure of the representation.\psp

\psp {\bf Case 2}. $r=\m$. \psp

In this case,
\begin{equation}
\Delta=-\sum\limits_{i=1}^{\m}x_i\partial_{x_{\m+i}}+\sum\limits_{j=1}^n\partial_{\theta_j}\partial_{\theta_{n+j}},\;\;
\eta=\sum\limits_{i=1}^{\m}x_{\m+i}\partial_{x_i}+\sum\limits_{j=1}^n\theta_j\theta_{n+j}.
\end{equation}
Furthermore, ${\cal H}^{\m}_k=0$ when $k\geq n$.

By similar arguments as those in Case 1, we can obtain that ${\cal
H}^{\m}_k$ is irreducible and ${\cal A}^{\m}_k={\cal
H}^{\m}_k\oplus\eta{\cal A}^{\m}_{k-2}$ when $k\leq n$. Set
\begin{eqnarray}
h(\vec{k},\vec{l},\vec{s})&=&\prod\limits_{t=1}^{\m}
x_t^{k_t}\prod\limits_{1\leq i<j\leq
\m}(x_ix_{\m+j}-x_jx_{\m+i})^{k_{i,j}}\prod\limits_{j=1}^n
\theta_j^{l_j}\theta_{n+j}^{l_{n+j}}\prod\limits_{1\leq i<j\leq
n}(\theta_i\theta_{n+i}\nonumber\\
&&-\theta_j\theta_{n+j})^{l_{i,j}}\prod\limits_{1\leq p\leq \m,1\leq
q\leq n}(x_{\m+p}-x_p\theta_q\theta_{n+q})^{s_{p,q}}
\end{eqnarray}
where
\begin{eqnarray}
&&\vec{k}=(k_1,\cdots,k_{\m};k_{1,2},k_{1,3},\cdots,k_{1,\m},k_{2,3}\cdots,k_{\m-1,\m})\in\mathbb{N}^{\frac{\m(\m+1)}{2}},\\
&&\vec{l}=(l_1,\cdots,l_n;l_{1,2},\cdots,l_{1,n},l_{2,3},\cdots,l_{2,n},\cdots,l_{n-1,n})\in\{0,1\}^{\frac{n(n+1)}{2}},\\
&&\vec{s}=(s_{1,1},\cdots,s_{1,n},\cdots,s_{\m,n})\in\{0,1\}^{\m n}.
\end{eqnarray}
Denote
\begin{eqnarray}
I&=&\{(\vec{k},\vec{l},\vec{s})\mid l_t+l_{n+t}+\sum\limits_{1\leq
i<t}l_{i,t}+\sum\limits_{t<j\leq
n}l_{t,j}+\sum\limits_{p=1}^{\m}s_{p,t}\leq 1\mbox{ for
}t\in\overline{1,n};\nonumber\\
&&k_{i,j}k_t=0\mbox{ for }i<j<t;
\; k_{i,j}k_{i',j'}=0\mbox{ for
}i>i'\mbox{ and }j<j';\nonumber\\
&&k_tl_{i,j}=0\mbox{ for }t\in\overline{1,\m},1\leq i<j\leq
n;
\;k_ts_{p,q}=0\mbox{ for }t<p;\;\nonumber\\
&& k_{i,j}s_{p,q}=0\mbox{ for }i<j<p;\;s_{p,q}s_{p',q'}=0\mbox{ for
}p>p'\mbox{ and }q<q';\nonumber\\
&& l_{i,j}=0 \mbox{ if
}l_t=l_{n+t}=\sum\limits_{p=1}^{\m}s_{p,t}=\sum\limits_{i'<t}l_{i',t}+\sum\limits_{j'>t}l_{t,j'}=0\mbox{
for some }i<t<j; \nonumber\\
&&l_{i,j}l_{i',j'}=0\mbox{ if }i<i'<j<j';\; l_{i,j}s_{p,q}=0\mbox{
if }i<j<q \}.
\end{eqnarray}
\begin{thm}
If $k\leq n$, the subspace ${\cal H}^{\m}_k$ is an irreducible
highest weight submodule. The highest weight is
$-2\lambda_{\m}+\nu_k$ (resp. $k\lambda_{\m-1}-(k+2)\lambda_{\m}$)
if $k>0$ (resp. $k\leq 0$) and a corresponding highest weight vector
is $\theta_1\cdots\theta_k$ (resp. $x_{\m}^{-k}$). Moreover, ${\cal
A}^{\m}_k={\cal H}^{\m}_k\oplus\eta{\cal A}^{\m}_{k-2}$ and
\begin{equation}
\begin{array}{l}
\{h(\vec{k},\vec{l},\vec{s})\mid (\vec{k},\vec{l},\vec{s})\in I;\;
\sum\limits_{t=1}^n(l_t+l_{n+t})+2\sum\limits_{1\leq i<j\leq
n}l_{i,j}+\sum\limits_{1\leq p\leq \m,1\leq q\leq
n}s_{p,q}-\sum\limits_{t=1}^{\m}k_t=k; \} \label{d7}
\end{array}
\end{equation}
forms a basis for ${\cal H}^{\m}_k$.
\end{thm}
\begin{proof} Denote by $V$ the subspace spanned by (\ref{d7}). It is easy to
check ${\cal H}^{\m}_k\supset V$. For the reverse inclusion, we
prove it by induct on $n$. Take any $0\neq f\in{\cal H}^{\m}_k$. We
write
\begin{equation}
f=f_0+f_1\theta_{n}+f_2\theta_{2n}+f_3\theta_{n}\theta_{2n}
\end{equation}
where
\begin{equation}f_i\in\mathbb{C}[x_1,\cdots,x_{2\m};\theta_1,\cdots,\theta_{n-1},\theta_{n+1}\cdots,\theta_{2n-1}].
\end{equation}
Let
\begin{equation}\Delta'=-\sum\limits_{i=1}^{\m}x_i\partial_{x_{\m+i}}+\sum\limits_{j=1}^{n-1}\partial_{\theta_j}\partial_{\theta_{n+j}}.
\end{equation}

Since
\begin{equation}
\Delta(f)=0=\Delta'(f_0)+\Delta'(f_1)\theta_n+\Delta'(f_2)\theta_n+\Delta'(f_3)\theta_n\theta_{2n}-f_3,
\end{equation}
we get
\begin{equation}\Delta'(f_0)-f_3=\Delta'(f_1)=\Delta'(f_2)=\Delta'(f_3)=0.\end{equation}
Thus by inductive assumption, we obtain
$f_1\theta_{n},f_2\theta_{2n}\in V$. We may assume
$f_3=h(\vec{k},\vec{l},\vec{s})$ with $(\vec{k},\vec{l},\vec{s})\in
I$ and
\begin{equation}
\sum\limits_{t=1}^n(l_t+l_{n+t})+2\sum\limits_{1\leq i<j\leq
n}l_{i,j}+\sum\limits_{1\leq p\leq \m,1\leq q\leq
n}s_{p,q}-\sum\limits_{t=1}^{\m}k_t=k-2,\end{equation}
\begin{equation}
l_n=l_{2n}=l_{1,n}=\cdots=l_{n-1,n}=s_{1,n}=\cdots=s_{\m,n}=0.
\end{equation}
Suppose that there exists some $t\in\overline{1,\m}$ such that
$k_t>0$. Let
\begin{equation}t_0=\mbox{min}\{t\mid k_t>0\},\;\;
f_0'=-x_{t_0}^{-1}x_{\m+t_0}f_3.\end{equation} We have
\begin{equation}
\Delta'(f_0')=f_3\end{equation} and
\begin{equation}0=\Delta(f_0+f_3\theta_n\theta_{2n})=\Delta'(f_0-f_0')+(\Delta'(f_0')-f_3).\end{equation}
So $f_0-f_0'\in V$ by inductive assumption and
\begin{equation}f_0'+f_3\theta_n\theta_{2n}=\Delta'(f_0)-f_3=-x_{t_0}^{-1}f_3(x_{\m+t_0}-x_{t_0}\theta_n\theta_{2n})\in
V.\end{equation} Thus $f\in V$.

Now we assume $k_t=0$ for all $t\in\overline{1,\m}$. Then  there
must exist $j\in\overline{1,n-1}$ such that
\begin{equation}l_j=l_{n+j}=\sum\limits_{1\leq i<j}l_{i,j}=\sum\limits_{j<t\leq
n}l_{j,t}=\sum\limits_{p=1}^{\m}s_{p,j}=0\end{equation} because
\begin{equation}
\sum\limits_{t=1}^{n-1}\big(l_t+l_{n+t}+\sum\limits_{1\leq
i<t}l_{i,t}+\sum\limits_{t<j\leq
n}l_{t,j}+\sum\limits_{p=1}^{\m}s_{p,t}\big)=k-2<n-1.
\end{equation}
Set
\begin{equation}j_0=\mbox{max}\{j\mid l_j=l_{n+j}=\sum\limits_{1\leq
i<j}l_{i,j}=\sum\limits_{j<t\leq
n}l_{j,t}=\sum\limits_{p=1}^{\m}s_{p,j}=0\}\end{equation}
 and
$f_0''=-f_3\theta_{j_0}\theta_{n+j_0}$. We have $\Delta'(f_0'')=f_3$
and
\begin{equation}
0=\Delta(f_0+f_3\theta_n\theta_{2n})=\Delta'(f_0-f_0')+\Delta'(f_0'')-f_3.
\end{equation}
So
\begin{equation}
f_0-f_0''\in
V,\;\;f_0''+f_3\theta_n\theta_{2n}=-f_3(\theta_{j_0}\theta_{n+j_0}-\theta_n\theta_{2n})\in
V,\end{equation} which implies $f\in V$.

Next we check the linear independence of (\ref{d7}). Again induct on
$n$. Suppose
\begin{equation}
\sum_{\stackrel{(\vec{k},\vec{l},\vec{s})\in
I}{\sum\limits_{t=1}^n(l_t+l_{n+t})+2\sum\limits_{1\leq i<j\leq
n}l_{i,j}+\sum\limits_{1\leq p\leq \m,1\leq q\leq
n}s_{p,q}-\sum\limits_{t=1}^{\m}k_t=k}}
a_{\vec{k},\vec{l},\vec{s}}h(\vec{k},\vec{l},\vec{s})=0.\label{d71}
\end{equation}
We write
\begin{eqnarray}
h(\vec{k},\vec{l},\vec{s})=h'(\vec{k},\vec{l},\vec{s})(\theta_i\theta_{n+i}-\theta_n\theta_{2n})\mbox{
if }l_{i,n}=1,\\
h(\vec{k},\vec{l},\vec{s})=h''(\vec{k},\vec{l},\vec{s})(x_{\m+p}-x_p\theta_n\theta_{2n})\mbox
{ if }s_{p,n}=1.
\end{eqnarray}
Thus (\ref{d71}) becomes
\begin{eqnarray}
&& \sum\limits_{(\vec{k},\vec{l},\vec{s})\in
I;l_n=1}a_{\vec{k},\vec{l},\vec{s}}h(\vec{k},\vec{l},\vec{s})
+\sum\limits_{(\vec{k},\vec{l},\vec{s})\in
I;l_{2n}=1}a_{\vec{k},\vec{l},\vec{s}}h(\vec{k},\vec{l},\vec{s})+\sum\limits_{\stackrel{(\vec{k},\vec{l},\vec{s})\in
I;}{l_n=l_{2n}=\sum l_{i,n}=\sum s_{p,n}=0}}
a_{\vec{k},\vec{l},\vec{s}}h(\vec{k},\vec{l},\vec{s})\nonumber\\
&&+\sum\limits_{i=1}^{n-1}\sum\limits_{(\vec{k},\vec{l},\vec{s})\in
I;l_{i,n}=1}
a_{\vec{k},\vec{l},\vec{s}}h'(\vec{k},\vec{l},\vec{s})(\theta_i\theta_{n+i}-\theta_n\theta_{2n})\nonumber\\
&&+\sum\limits_{p=1}^{\m}\sum\limits_{(\vec{k},\vec{l},\vec{s})\in
I;s_{p,n}=1}
a_{\vec{k},\vec{l},\vec{s}}h''(\vec{k},\vec{l},\vec{s})(x_{\m+p}-x_p\theta_n\theta_{2n})=0.
\end{eqnarray}
We get
\begin{equation}
a_{\vec{k},\vec{l},\vec{s}}=0\mbox{ if }l_n=l_{2n}=\sum l_{i,n}=\sum
s_{p,n}=0 \mbox{ or }l_n=1\mbox{ or }l_{2n}=1,
\end{equation}
and
\begin{equation}
\sum\limits_{i=1}^{n-1}\sum\limits_{(\vec{k},\vec{l},\vec{s})\in
I;l_{i,n}=1}
a_{\vec{k},\vec{l},\vec{s}}h'(\vec{k},\vec{l},\vec{s})+\sum\limits_{p=1}^{\m}\sum\limits_{(\vec{k},\vec{l},\vec{s})\in
I;s_{p,n}=1}
a_{\vec{k},\vec{l},\vec{s}}x_ph''(\vec{k},\vec{l},\vec{s})=0.
\end{equation}
Since $h'(\vec{k},\vec{l},\vec{s})$ and
$x_ph''(\vec{k},\vec{l},\vec{s})$ are linearly independent by
inductive assumption, we get $a_{\vec{k},\vec{l},\vec{s}}=0$ for
$l_{i,n}=1$ or $s_{p,n}=1$.
\end{proof}
\psp

\textbf{Case 3.} $0<r<\m$ \psp

This case is a little more complicated. To study the structure of
the submodules ${\cal A}^r_k$ and ${\cal H}^r_k$, we first introduce
some subalgebras of $osp(2\m,2n)$. Set
\begin{equation}
\begin{array}{l}
L_1=\sum\limits_{r+1\leq i,j\leq
\m}\mathbb{C}(E_{i,j}-E_{\m+j,\m+i})+\sum\limits_{r+1\leq i<j\leq
\m}\big(\mathbb{C}(E_{i,\m+j}-E_{j,\m+i})\\
+\mathbb{C}(E_{\m+j,i}-E_{\m+i,j})\big)
\end{array}
\end{equation}
and
\begin{equation}
\begin{array}{l}
L'_1=\sum\limits_{r+1\leq i,j\leq
\m-1}\mathbb{C}(E_{i,j}-E_{\m+j,\m+i})+\sum\limits_{r+1\leq i<j\leq
\m-1}\big(\mathbb{C}(E_{i,\m+j}-E_{j,\m+i})\\
+\mathbb{C}(E_{\m+j,i}-E_{\m+i,j})\big).
\end{array}
\end{equation}
Denote $L_1^+=osp(2\m,2n)^+\cap L_1$ and ${L'}_1^+=osp(2\m,2n)^+\cap
L'_1$. We treat $L'_1=0$ when $r=\m-1$. Let
\begin{eqnarray}
L_2&=&\sum\limits_{i,j=1}^r\big(\mathbb{C}(E_{i,j}-E_{\m+j,\m+i})+\mathbb{C}(E_{i,\m+j}-E_{j,\m+i})\nonumber\\
&&+\mathbb{C}(E_{\m+i,j}-E_{\m+j,i})\big)
+\sum\limits_{p,q=1}^n\big(\mathbb{C}(E_{2\m+p,2\m+q}-E_{2\m+n+q,2\m+n+p})\nonumber\\
&&+\mathbb{C}(E_{2\m+p,2\m+n+q}+E_{2\m+q,2\m+n+p})+\mathbb{C}(E_{2\m+n+p,2\m+q}+E_{2\m+n+q,2\m+p})\big)\nonumber\\
&&+
\sum\limits_{i\in\overline{1,r};p\in\overline{1,n}}\big(\mathbb{C}(E_{i,2\m+p}-E_{2\m+n+p,\m+i})+\mathbb{C}(E_{i,2\m+n+p}+E_{2\m+p,\m+i})\nonumber\\
&&+\mathbb{C}(E_{\m+i,2\m+p}-E_{2\m+n+p,i})+\mathbb{C}(E_{\m+i,2\m+n+p}+E_{2\m+p,i})
\end{eqnarray}
and
\begin{eqnarray}
L_2^+&=&L_2\cap osp(2\m,2n)_{0}^++\sum\limits_{1\leq i\leq m,1\leq
p\leq
n}\big(\mathbb{C}(E_{\m+i,2\m+n+p}+E_{2\m+p,i})\nonumber\\
&&+\mathbb{C}(E_{i,2\m+n+p}+E_{2\m+p,\m+i})\big).
\end{eqnarray}

We have the following result:
\begin{thm}
When $0<r<\m$ and $k\leq n-\m+r+1$, the submodule ${\cal H}^r_k$ is
irreducible and ${\cal A}^r_k={\cal H}^r_k\bigoplus\eta{\cal
A}^r_{k-2}$. If $k>n-\m+r+1$, we have the following composition
series
\begin{eqnarray}
&&{\cal H}^r_k\supset\eta^{k-n+\m-r-1}{\cal
H}^r_{-k+2(n-\m+r+1)}\supset\{0\}\mbox{ if }r<\m-1;\label{d6252}\\
&&{\cal H}^{\m-1}_k\supset\la x_{\m}^k\ra\supset\eta^{k-n}{\cal
H}^{\m-1}_{-k+2n}\supset\{0\}.\label{d6251}
\end{eqnarray}
The subspace ${\cal H}^r_k$ ($k\in\mathbb{Z}$) has a basis
\begin{equation}
\begin{array}{l}
\{\sum\limits_{\stackrel{r_1,\cdots,r_{\m-1}\in\mathbb{N};}{l_1,\cdots,l_n\in\{0,1\}}}\frac{(-1)^{\sum\limits_{i=r+1}^{\m-1}r_i+\sum\limits_{j=1}^nl_j}\prod\limits_{i=1}^r
\left(
\begin{array}{l}
\alpha_{\m+i}\\
r_i
\end{array}
\right)\prod\limits_{i=r+1}^{\m}\left(
\begin{array}{l}
\alpha_i\\
r_i
\end{array}
\right)\left(
\begin{array}{l}
\alpha_{\m+i}\\
r_i
\end{array}
\right)\prod\limits_{j=1}^n \left(
\begin{array}{c}
\beta_j\\
l_j
\end{array}
\right)\left(
\begin{array}{c}
\beta_{n+j}\\
l_j
\end{array}
\right)} {\left(
\begin{array}{c}
\alpha_{\m}+\sum\limits_{i=1}^{\m-1}r_i+\sum\limits_{j=1}^nl_j\\
\alpha_{\m}
\end{array}
\right)\left(
\begin{array}{c}
\alpha_{2\m}+\sum\limits_{i=1}^{\m-1}r_i+\sum\limits_{j=1}^nl_j\\
\alpha_{2\m}
\end{array}
\right)\left(
\begin{array}{c}
\sum\limits_{i=1}^{\m-1}r_i+\sum\limits_{j=1}^nl_j\\
r_1,\cdots,r_{\m-1}
\end{array}
\right)\prod\limits_{i=1}^rr_i!}\\
\times\prod\limits_{i=1}^{r}x_i^{\alpha_i+r_i}x_{\m+i}^{\alpha_{\m+i}-r_i}\prod\limits_{i=r+1}^{\m}x_i^{\alpha_i-r_i}x_{\m+i}^{\alpha_{\m+i}-r_i}\prod\limits_{j=1}^n(-1)^{\beta_jl_j}\theta_j^{\beta_j-l_j}\theta_{n+j}^{\beta_{n+j}-l_j}x_{\m}^{\alpha_\m+\sum\limits_{i=1}^{\m-1}r_i+\sum\limits_{j=1}^nl_j}\\
\times
x_{2\m}^{\alpha_{2\m}+\sum\limits_{i=1}^{\m-1}r_i+\sum\limits_{j=1}^nl_j}
\mid \alpha_i\in\mathbb{N},\beta_j\in\{0,1\};
\sum\limits_{j=1}^{2n}\beta_j-\sum\limits_{i=1}^{r}\alpha_i+\sum\limits_{i=r+1}^{2\m}\alpha_i=k;\alpha_{\m}\alpha_{2\m}=0\}.
\end{array}
\end{equation}
\end{thm}
\begin{proof}
\textit{(1) The subspace ${\cal H}^{\m-1}_k$ is generated by
\begin{equation}
f_{l,p,s}=\sum\limits_{i=0}^l\frac{l!(l+p)!}{i!(l-i)!(l+p-i)!}x_{\m}^p(x_{\m}x_{2\m})^{l-i}x_{\m-1}^{s-i}x_{2\m-1}^i\theta_1\cdots\theta_n\label{d6231}
\end{equation}
and
\begin{equation}
g_{l,p,s}=\sum\limits_{i=0}^l\frac{l!(l+p)!}{i!(l-i)!(l+p-i)!}x_{2\m}^p(x_{\m}x_{2\m})^{l-i}x_{\m-1}^{s-i}x_{2\m-1}^i\theta_1\cdots\theta_n\label{d6232}
\end{equation}
as an $(L_1+L_2)$-submodule with $l,p,s\in\mbb{N}$, $p+2l-s=k$ and
$l\leq s$.}\psp

Using Lemma 2.1, we obtain that the subspace ${\cal H}^{\m-1}_k$ is
spanned by
\begin{equation}
\sum\limits_{i=0}^\infty(-\int_{(\m)}\int_{(2\m)})^ix_{\m}^{\alpha_{\\m}}x_{2\m}^{\alpha_{2\m}}(-\sum\limits_{j=1}^{\m-1}x_j\partial_{x_{\m+j}}+\sum\limits_{j=1}^n\partial_{\theta_j}\partial_{\theta_{n+j}})^i(g),
\end{equation}
where
\begin{equation}\alpha_{\m}\alpha_{2\m}=0,\qquad
g\in\mathbb{C}[x_1,\cdots,x_{\m-1},x_{\m+1},\cdots,x_{2\m-1};\theta_1,\cdots,\theta_{2n}]\end{equation}
and
\begin{equation}\int_{(t)}x^\alpha=\frac{x_tx^\alpha}{\alpha_t+1}\ \
\mbox{for}\;\;t=\m,2\m.\end{equation}

According to Theorem 2.3, we can write $g=\sum_l
X_l(x_{\m-1}^{s-l}x_{2\m-1}^l\theta_1\cdots\theta_n)$, where $X_l\in
U(L_2)$, $2l-s+\alpha_\m+\alpha_{2\m}+n=k$. Thus ${\cal H}^{\m-1}_k$
is spanned by
\begin{equation}
X_l\sum\limits_{i=0}^l\frac{p!l!}{i!(l-i)!(p+l-i)!}x_{\m}^p(x_{\m}x_{2\m})^{l-i}x_{\m-1}^{s-i}x_{2\m-1}^i\theta_1\cdots\theta_n=X_l\frac{p!}{(l+p)!}f_{l,p,s}
\end{equation}
and
\begin{equation}
X_l\sum\limits_{i=0}^l\frac{p!l!}{i!(l-i)!(p+l-i)!}x_{2\m}^p(x_{\m}x_{2\m})^{l-i}x_{\m-1}^{s-i}x_{2\m-1}^i\theta_1\cdots\theta_n=X_l\frac{p!}{(l+p)!}g_{l,p,s},
\end{equation}
where $2l+p+n-s=k$. Consequently, ${\cal H}^{\m-1}_k$ is generated
by (\ref{d6231}) and (\ref{d6232}) as an $L_2$-submodule.

\textit{(2) As an $(L_1+L_2)$-submodule, ${\cal H}^r_k$ ($r<\m-1$)
is generated by
\begin{equation}
h_{l,p,s}=\sum\limits_{i=0}^l\frac{l!(l+p+\m-r-1)!}{i!(l-i)!(l+p+\m-r-1)!}x_{r+1}^p
(\sum\limits_{j=r+1}^{\m}x_jx_{\m+j})^{l-i}x_r^{s-i}x_{\m+r}^i\theta_1\cdots\theta_n
\end{equation}
for $2l+p+n-s=k$ and $l\leq s$.}

Again by Lemma 2.1, we obtain that ${\cal H}^r_k$ ($r<\m-1$) is
spanned by
\begin{eqnarray}
&&X_{l,p,s,t}\sum\limits_{i=0}^l(-\int_{(\m)}\int_{(2\m)})^ix_{\m}^{\alpha_{\m}}x_{2\m}^{q-\alpha_{\m}}(-x_r\partial_{x_{\m+r}}+\sum\limits_{j=r+1}^{\m-1}\partial_{x_j}\partial_{x_{\m+j}})^i(x_{r+1}^{\alpha_{r+1}}x_{\m+r+1}^{p-\alpha_{r+1}}\nonumber\\
&&\times u'^{l-t}x_r^{s-t}x_{\m+r}^t\theta_1\cdots\theta_n),
\end{eqnarray}
where
\begin{equation}u'=\sum\limits_{j=r+1}^{\m-1}x_jx_{\m+j},\;\;\alpha_{\m}\in\{0,q\},\;\; \alpha_{r+1}\in\{0,p\},\;\;2l+p+q-s+n=k\end{equation}
 and
$X_{l,p,s,t}\in U(L'_1+L_2)$. So ${\cal H}^r_k$ is generated by
\begin{eqnarray}
g_{p,q,s,t,\alpha_{r+1},\alpha_{\m}}&=&\sum\limits_{i=0}^l(-\int_{(\m)}\int_{(2\m)})^ix_{\m}^{\alpha_{\m}}x_{2\m}^{q-\alpha_{\m}}(-x_r\partial_{x_{\m+r}}+\sum\limits_{j=r+1}^{\m-1}\partial_{x_j}\partial_{x_{\m+j}})^i\nonumber\\
&&(x_{r+1}^{\alpha_{r+1}}x_{\m+r+1}^{p-\alpha_{r+1}}u'^{l-t}x_r^{s-t}x_{\m+r}^t\theta_1\cdots\theta_n)
\end{eqnarray}
as an $(L'_1+L_2)$-submodule. Denote by
$\lambda_{p,q,s,\alpha_{r+1},\alpha_{\m}}$ the weight of
$g_{p,q,s,t,\alpha_{r+1},\alpha_{\m}}$. Note
\begin{equation}
{L'}_1^+(g_{p,q,s,t,\alpha_{r+1},\alpha_{\m}})=L_2^+(g_{p,q,s,t,\alpha_{r+1},\alpha_{\m}})=0.
\end{equation}
Hence
\begin{equation}
{\cal H}^r_k=\sum U(L'_1+L_2)(g_{p,q,s,t,\alpha_{r+1},\alpha_{\m}})
\end{equation}
and
\begin{equation}
\begin{array}{ll}
&({\cal H}^r_k)_{\lambda_{p,q,s,\alpha_{r+1},\alpha_{\m}}}\cap
\mbox{Span}\{g\theta_1\cdots\theta_n\mid g\in\mathbb{C}[x_{r+1},x_{\m+r+1},u',x_{\m},x_{2\m}]\}\\
=&\mbox{Span}\{g_{p,q,s,t,\alpha_{r+1},\alpha_{\m}}\mid 0\leq t\leq
\mbox{min}\{s,l\}\}.
\end{array}
\end{equation}
Note
\begin{equation}
\mbox{dim}\;\mbox{Span}\;\{g_{l,p,s,t,\alpha_{r+1},\alpha_{\m}}\mid
0\leq t\leq \mbox{min}\{s,l\}\}=\left\{
\begin{array}{l}
l+1\mbox{ if }l\leq s,\\
s+1\mbox{ if }l>s.
\end{array}
\right.
\end{equation}
On the other hand,
\begin{eqnarray}
&&(\sum_{t=r+1}^{\m}\partial_{x_t}\partial_{x_{\m+t}})\big(\sum_{j=0}^d\frac{(-1)^j}{j!(d-j)!(d+p+\m-r-2-j)!(q+j)!}x_{r+1}^{\alpha_{r+1}}x_{\m+r+1}^{p-\alpha_{r+1}}\nonumber\\
&&\times
x_{\m}^{\alpha_{\m}}x_{2\m}^{q-\alpha_{\m}}u'^{d-j}(x_{\m}x_{2\m})^j\big)=0,
\end{eqnarray}
which means
\begin{eqnarray}
&&\sum_{j=0}^d\frac{(-1)^j}{j!(d-j)!(d+p+\m-r-2-j)!(q+j)!}x_{r+1}^{\alpha_{r+1}}x_{\m+r+1}^{p-\alpha_{r+1}}\nonumber\\
&&\times
x_{\m}^{\alpha_{\m}}x_{2\m}^{q-\alpha_{\m}}u'^{d-j}(x_{\m}x_{2\m})^j\in
U(L_1)(x_{r+1}^{p+q+2d}).
\end{eqnarray}
So we get
\begin{eqnarray}
&&\sum\limits_{i=0}^{l-d}\frac{l!(l+d+p+\m-r-1)!}{i!(l-d-i)!(l+d+p+\m-r-1)!}\big(\sum\limits_{j=0}^d\frac{(-1)^j}{j!(d-j)!(q+j)!(p+d+\m-r-2-j)!}\nonumber\\
&&\times
u'^{d-j}(x_{\m}x_{2\m})^j\big)x_{r+1}^{\alpha_{r+1}}x_{\m+r+1}^{p-\alpha_{r+1}}x_{\m}^{\alpha_{\m}}x_{2\m}^{q-\alpha_{\m}}(u'+x_{\m}x_{2\m})^{l-d-i}x_r^{s-i}x_{\m+r}^i\theta_1\cdots\theta_n\nonumber\\
&&\in U(L_1)(h_{l-d,p+q+2d,s})\cap
\mbox{Span}\;\{g_{l,p,s,t,\alpha_{r+1},\alpha_{\m}}\mid 0\leq t\leq
\mbox{min}\{s,l\}\},
\end{eqnarray}
where $0\leq d\leq \mbox{min}\{l,s\}$. Thus
\begin{eqnarray}
&&\mbox{dim}\big(\bigoplus_{d=0}^{\mbox{\small min}\{l,s\}}
U(L_1)(h_{l-d,p+q+2d,s})\cap
\mbox{Span}\;\{g_{l,p,s,t,\alpha_{r+1},\alpha_{\m}}\mid 0\leq t\leq
\mbox{min}\{s,l\}\}\nonumber\\
=&&\left\{
\begin{array}{l}
l+1\mbox{ if }l\leq s,\\
s+1\mbox{ if }l>s.
\end{array}
\right.
\end{eqnarray}
Therefore, we have
\begin{eqnarray}
&&\bigoplus\limits_{d=0}^{\mbox{\small min}\{l,s\}}
U(L_1)(h_{l-d,p+q+2d,s})\cap\mbox{Span}\;\{g_{l,p,s,t,\alpha_{r+1},\alpha_{\m}}\mid
0\leq t\leq
\mbox{min}\{s,l\}\}\nonumber\\
=&&\mbox{Span}\;\{g_{l,p,s,t,\alpha_{r+1},\alpha_{\m}}\mid 0\leq
t\leq \mbox{min}\{s,l\}\},
\end{eqnarray}
which implies
\begin{equation}
g_{l,p,s,t,\alpha_{r+1},\alpha_{\m}}\in\bigoplus_{d=0}^{\mbox{\small
min}\{l,s\}} U(L_1)(h_{l-d,p+q+2d,s}).\end{equation}
 Hence we are done.\psp

\textit{(3) We claim that
\begin{equation}
{\cal H}^r_k=\left\{
\begin{array}{lll}
\la x_r^{-k}\ra  &\mbox{\it if} &k\leq0,\\
\la x_{r+1}^k\ra &\mbox{\it if} &k>0,r<\m-1,\\
\la x_{\m}^k\ra+\la x_{2\m}^k\ra &\mbox{\it if} &k>0,r=\m-1.
\end{array}
\right.\label{d31}
\end{equation}}
\psp

When $r=\m-1$, we have
\begin{eqnarray}
f_{l,p,s}&=&\sum\limits_{i=0}^l\frac{l!(l+p)!}{i!(l-i)!(l+p-i)!}x_{\m}^p(x_{\m}x_{2\m})^{l-i}x_{\m-1}^{s-i}x_{2\m-1}^i\theta_1\cdots\theta_n\nonumber\\
&=&(E_{2\m-1,\m}-E_{2\m,\m-1})^l(x_{\m-1}^{s-l}x_{\m}^{p+l}\theta_1\cdots\theta_n)
\end{eqnarray}
and
\begin{eqnarray}
g_{l,p,s}&=&\sum\limits_{i=0}^l\frac{l!(l+p)!}{i!(l-i)!(l+p-i)!}x_{2\m}^p(x_{\m}x_{2\m})^{l-i}x_{\m-1}^{s-i}x_{2\m-1}^i\theta_1\cdots\theta_n\nonumber\\
&=&(E_{2\m-1,2\m}-E_{\m,\m-1})^l(x_{\m-1}^{s-l}x_{2\m}^{p+l}\theta_1\cdots\theta_n).
\end{eqnarray}
It is straightforward to check
\begin{equation}
x_{\m-1}^{s-l}x_{\m}^{p+l}\theta_1\cdots\theta_n,\;x_{\m-1}^{s-l}x_{2\m}^{p+l}\theta_1\cdots\theta_n\in\left\{
\begin{array}{ll}
\la x_{\m}^k\ra+\la x_{2\m}^k\ra&\mbox{if }k>0,\\
\la x_{\m-1}^{-k}\ra&\mbox{if }k\leq 0.
\end{array}
\right.
\end{equation}
Now we assume $r<\m-1$. Then $h_{0,p,p+n-k}\in\la x_{r+1}^k\ra$ (or
$\la x_r^{-k}\ra$). Note
\begin{equation}
\begin{array}{ll}
l(p+1)h_{l,p,s}=&\sum\limits_{j=r+2}^{\m}(E_{\m+j,r+1}-E_{\m+r+1,j})(-E_{j,r}+E_{\m+r,\m+j})(h_{l-1,s-1,p+1})\\
&+(\m-r)(E_{\m+r,r+1}-E_{\m+r+1,r})(h_{l-1,p+1,s-1}).
\end{array}\end{equation}
So $h_{l,p,s}\in\la x_{r+1}^k\ra$ (or $\la x_r^{-k}\ra$) by
inductive assumption. Therefore, (\ref{d31}) holds. \psp

\textit{(3) When $k\leq n-\m+r+1$, the submodule ${\cal H}^r_k$ is
irreducible and ${\cal A}^r_k={\cal H}^r_k\oplus\eta{\cal
A}^r_{k-2}$.}

\psp We may assume $r<\m-1$ and $k>0$. The proof for $r=\m-1$ or
$k\leq0$ is similar. For any submodule $W\subset{\cal H}^r_k$, there
should be some weight vector $g\in W$ such that $L_1^+(g)=0$ and
$L_2^+(g)=0$. Thus
$g=\sum\limits_{i=0}^la_ix_{r+1}^pu^{l-i}x_r^{s-i}x_{\m+r}^i\theta_1\cdots\theta_n$
for some $a_i\in\mathbb{C}$. Since
\begin{eqnarray}
0&=&\Delta(g)=-\sum_{i=1}^la_iix_{r+1}^pu^{l-i}x_r^{s-i+1}x_{\m+r}^{i-1}\theta_1\cdots\theta_n\nonumber\\
&&+\sum\limits_{i=0}^{l-1}x_{r+1}^pu^{l-i}x_r^{s-i}x_{\m+r}^i\theta_1\cdots\theta_n
\end{eqnarray}
we get
\begin{equation}
a_{i+1}(i+1)=a_i(l-i)(p+l+\m-r-i-1).
\end{equation}
Thus $g$ is a scalar multiple of
\begin{equation}
h_{l,p,s}=\sum\limits_{i=0}^l\frac{l!(l+p+\m-r-1)!}{i!(l-i)!(l+p+\m-r-1-i)!}x_{r+1}^pu^{l-i}x_r^{s-i}x_{\m+r}^i\theta_1\cdots\theta_n.
\end{equation}
If $l>0$, we have
\begin{equation}
(E_{r,\m+r+1}-E_{r+1,\m+r})(h_{l,p,s})=-l(l+p+\m-r-1-s)h_{l-1,p+1,s-1},
\end{equation}
and
\begin{equation}
l+\m-r-1+p-s=k-n-l+\m-r-1\leq -l<0,\end{equation} which implies
$h_{l-1,p+1,s-1},\cdots,h_{0,p+l,s-l}\in W$. It is easy to see
$x_{r+1}^k\in\la h_{0,p+l,s-l}\ra$. Hence $W={\cal H}^r_k$.

By the similar arguments as those in  (3) of the proof of Theorem
2.2, we obtain ${\cal A}^r_k={\cal H}^r_k\oplus\eta{\cal
A}^r_{k-2}$.\psp

\textit{(4)(\ref{d6252}) and (\ref{d6251}) are composition series.}

\psp According to (2),
\begin{equation}{\cal H}^{\m-1}_k=\bigoplus_{l,p,s=0,p+2l-s=k,l\leq s}^\infty
(U(L_2)f_{l,p,s}\oplus U(L_2)g_{l,p,s})\end{equation} and
\begin{equation}{\cal H}^r_k= \bigoplus_{l,p,s=0,p+2l-s=k,l\leq s}^\infty U(L_1+L_2)h_{l,p,s}.\end{equation} Again by
 the similar arguments as those in (4) of the proof of Theorem 2.2, we can get the composition
series.
\end{proof}

\section{Proof of Theorem 2}
\setcounter{equation}{0} In this section, we  discuss the
$osp(2\m,2n)$-module ${\cal A}'_k$ defined in (\ref{I6181}) and
(\ref{i61}). The following facts will be used.
\begin{prop}
If $n=0$, then the subspace ${\cal A}'_k$ ($k\neq \m$) is an
irreducible $so(2\m,\mathbb{C})$-submodule and ${\cal A}'_{\m}=\la
\theta_1\cdots\theta_{\m}\ra\oplus\la
\theta_1\cdots\theta_{\m-1}\theta_{2\m}\ra$.

If $\m=0$, then the subspace ${\cal A}'_k$ is an irreducible
$sp(2n,\mathbb{C})$-submodule when $S_1\cup T_1\neq\emptyset$ or
$k\neq0$. When $S_1\cup T_1=\emptyset$, we can assume $T=\ol{1,n}$
by symmetry. In this case, ${\cal A}'_0=\la 1\ra\oplus\la
x_{n-1}x_{2n}-x_nx_{2n-1}\ra$. (cf. [\ref{r6}])
\end{prop}\psp

Now let us deal with the general case with $\m>0$ and $n>0$.

In fact, if $S_1\neq\emptyset$, we can take a $j_0\in S_1$ and
$0\neq f\in{\cal A}'_k$. Since
\begin{equation}
E_{2\m+j_0,2\m+n+j_0}\mid_{{\cal A}'}=x_{j_0}\partial_{x_{n+j_0}},
\end{equation}
we can assume $\partial_{x_{n+j_0}}(f)=0$. Applying
\begin{equation}
(E_{\m+i,2\m+n+j_0}+E_{2\m+j_0,i})\mid_{{\cal
A}'}=\theta_{\m+i}\partial_{x_{n+j_0}}+x_{j_0}\partial_{\theta_i}
\end{equation}
and
\begin{equation}
(E_{i,2\m+n+j_0}+E_{2\m+j_0,\m+i})\mid_{{\cal
A}'}=\theta_i\partial_{x_{n+j_0}}+x_{j_0}\partial_{\theta_{\m+i}}
\end{equation}
($i\in\overline{1,\m}$), we get a nonzero element $
f'=f'(x_1,\cdots,\hat{x_{n+j_0}},\cdots,x_{2n})\in\la f\ra$. Since
\begin{equation}
\mbox{Span}\;\{x^\alpha\mid x^\alpha\in{\cal A}'_k\}\end{equation}
is an irreducible $sp(2n,\mathbb{C})$-submodule according to
Proposition 3.1, we obtain
\begin{equation}
\mbox{Span}\;\{x^\alpha\mid x^\alpha\in{\cal A}'_k\}\subset\la f\ra.
\end{equation}
Observe
\begin{equation}
(E_{i,2\m+j_0}-E_{2\m+n+j_0,\m+i})|_{{\cal
A}'}=\theta_i\partial_{x_{j_0}}-x_{n+j_0}\partial_{\theta_{\m+i}}
\end{equation}
and
\begin{equation}
(E_{\m+i,2\m+j_0}-E_{2\m+n+j_0,i})|_{{\cal
A}'}=\theta_{\m+i}\partial_{x_{j_0}}-x_{n+j}\partial_{\theta_i}.
\end{equation}
Thus by induction on $t$, we can obtain
$x^\alpha\theta_{i_1}\cdots\theta_{i_t}\in\la f\ra$ for all
$i_1,\cdots,i_t\in\overline{1,2\m}$ and $\sum\limits_{i\in\bar{
T}}\alpha_i-\sum\limits_{j\in T}\alpha_j=k-t$. So $\la f\ra={\cal
A}'_k$, which implies that ${\cal A}'_k$ is irreducible. It can be
similarly proved when $T_1\neq\emptyset$.

\begin{thm}
1) The submodule ${\cal A}'_k$ is irreducible when $S_1\cup
T_1\neq\emptyset$. In particular, ${\cal A}'_k$ is not highest
weight type if $S_1\neq\emptyset$ and $T_1\neq\emptyset$.\\
2) If $S_1=\emptyset$ and $T_1=\emptyset$, we may assume
$T=\overline{1,n}$ by symmetry.\\
a) The submodule ${\cal A}'_k$ is irreducible and of highest weight
type when $k\neq
\m$. A highest weight vector is $x_n^{\m-k}\theta_1\cdots\theta_m$ (resp. $x_{2n}^{k-m_1}\theta_1\cdots\theta_m$) if $k>m_1$ (resp. $k<m_1$).\\
b) The submodule ${\cal A}'_{\m}=\la
\theta_1\cdots\theta_{\m}\ra\oplus\la
(x_{n-1}x_{2n}-x_nx_{2n-1})\theta_1\cdots\theta_{\m}\ra$ is a sum of
two irreducible submodules.
\end{thm}
\begin{proof}
Assume $T=\overline{1,n}$.

a) We claim that ${\cal A}'_k=\la
x_n^{\m-k}\theta_1\cdots\theta_{\m}\ra$ when $k<\m$.

In fact, we have
\begin{equation}
x_n^{l-k}\theta_1\cdots\theta_l=\frac{(l-k)!}{(\m-k)!}\prod_{t=l+1}^{\m}(-1)^{t-1}(E_{\m+t,2\m+2n}+E_{2\m+n,t})(x_n^{\m-k}\theta_1\cdots\theta_{\m})
\end{equation}
for $l\in\overline{1,\m}$. Thus we get
\begin{equation}
\mbox{Span}\;\{x^\alpha\theta_{i_1}\cdots\theta_{i_l}\mid
i_1,\cdots,i_l\in\overline{1,2m};\;
l+\sum_{i\in\bar{T}}\alpha_i-\sum_{i\in T}\alpha_i=k\}\subset\la
x_n^{\m-k}\theta_1\cdots\theta_{\m}\ra
\end{equation}
for $k<l<\m$ and $\theta_{i_1}\cdots\theta_{i_k}\in\la
x_n^{\m-k}\theta_1\cdots\theta_{\m}\ra$ by applying
$so(2\m,,\mathbb{C})$ and $sp(2n,\mathbb{C})$ to
$x_n^{l-k}\theta_1\cdots\theta_l$ (We treat
$\theta_1\cdots\theta_k=0$ if $k\leq0$). Since
\begin{eqnarray}
&&(E_{2\m,2\m+2n}+E_{2\m+n,\m})(-E_{2\m+2n,2\m+n})(x_n^{\m-k}\theta_1\cdots\theta_{\m})\nonumber\\
=&&x_n^{\m-k+1}\theta_{2\m}\theta_1\cdots\theta_{\m}+(-1)^{\m-1}(\m-k+1)x_n^{\m-k}x_{2n}\theta_1\cdots\theta_{\m-1},
\end{eqnarray}
we get $x_n^{\m-k+1}\theta_1\cdots\theta_{\m}\theta_{2\m}\in\la
x_n^{\m-1}\theta_1\cdots\theta_{\m}\ra$. Now we have
\begin{eqnarray}& &
x_n^{\m-k}\theta_1\cdots\theta_{\m-1}\theta_{2\m}\nonumber\\
&=&\frac{(-1)^{\m-1}}{\m+1-k}(E_{2\m,2\m+2n}+E_{2\m+n,\m})(x_n^{\m+1-k}\theta_1\cdots\theta_{\m}\theta_{2\m}).
\end{eqnarray}
Applying
\begin{equation}
(E_{\m+t,2\m+n}-E_{2\m+2n,t})|_{{\cal
A}'}=-x_n\theta_{\m+t}-x_{2n}\partial_{\theta_t}\label{62}
\end{equation}
and taking induction on $l$, we obtain
\begin{equation}
\mbox{Span}\;\{x^\alpha\theta_{i_1}\cdots\theta_{i_l}\mid
i_1,\cdots,i_l\in\overline{1,2m};\;
l+\sum_{i\in\bar{T}}\alpha_i-\sum_{i\in T}\alpha_i=k\}\subset\la
x_n^{\m-k}\theta_1\cdots\theta_{\m}\ra
\end{equation}
for $l\geq\m$. Since
\begin{eqnarray}
&&(x_{n-1}x_{2n}-x_nx_{2n-1})\theta_1\cdots\theta_k\nonumber\\
&=&(-1)^k(E_{\m+k+1,2\m+n-1}-E_{2\m+2n-1,k+1})(x_n\theta_1\cdots\theta_{k+1})\nonumber\\
&&-(-1)^k(E_{\m+k+1,2\m+n}-E_{2\m+2n,k+1})(x_{n-1}\theta_1\cdots\theta_{k+1}),
\end{eqnarray}
we get
\begin{equation}
\mbox{Span}\;\{x^\alpha\theta_{i_1}\cdots\theta_{i_l}\mid
i_1,\cdots,i_l\in\overline{1,2m};\;l+\sum_{i\in\bar{T}}\alpha_i-\sum_{i\in
T}\alpha_i=k\}\subset\la x_n^{\m-k}\theta_1\cdots\theta_{\m}\ra
\end{equation}
for $l=k$. Now by induction on $k-l$ for $l\leq k$ and  (\ref{62}),
we attain
\begin{equation}
\mbox{Span}\;\{x^\alpha\theta_{i_1}\cdots\theta_{i_l}\mid
i_1,\cdots,i_l\in\overline{1,2m};\;l+\sum_{i\in\bar{T}}\alpha_i-\sum_{i\in
T}\alpha_i=k\}\subset\la x_n^{\m-k}\theta_1\cdots\theta_{\m}\ra
\end{equation}
for $0\leq l<k$. Hence ${\cal A}'_k=\la
x_n^{\m-k}\theta_1\cdots\theta_{\m}\ra$.

Note that all the weight vectors annihilated by $osp(2\m,2n)^+_0$
are scalar multiples of
$x_n^{\m-k}\theta_1\cdots\theta_{\m-1}\theta_{2\m}$, $
x_n^{i-k}\theta_1\cdots\theta_i\;(k\leq i\leq \m)$,
$x_{2n}^{k-l}\theta_1\cdots\theta_l\;(0\leq l<k)$ and
$(x_{n-1}x_{2n}-x_nx_{2n-1})\theta_1\cdots\theta_k$. Since
\begin{equation}
\prod_{t=l+1}^k(E_{t,2\m+2n}+E_{2\m+n,\m+t})(x_{2n}^{k-l}\theta_1\cdots\theta_l)=(-1)^{l(k-l)}(k-l)!\theta_1\cdots\theta_k,
\end{equation}
\begin{eqnarray}& &
(E_{\m,2\m+2n}+E_{2\m+n,2\m})(x_n^{\m-k}\theta_1\cdots\theta_{\m-1}\theta_{2\m})\nonumber\\
& &=(-1)^{\m-1}(\m-k)x_n^{\m-k-1}\theta_1\cdots\theta_{\m-1},
\end{eqnarray}
\begin{eqnarray}
&&(E_{k+1,2\m+2n-1}+E_{2\m+n-1,\m+k+1})\big((x_{n-1}x_{2n}-x_nx_{2n-1})\theta_1\cdots\theta_k\big)\nonumber\\
&=&(-1)^{k+1}x_n\theta_1\cdots\theta_{k+1}
\end{eqnarray}
and
\begin{eqnarray}
\prod\limits_{j=i+1}^{\m}(-1)^j(E_{j,2\m+n}-E_{2\m+2n,\m+j})(x_n^{i-k}\theta_1\cdots\theta_i)=
x_n^{\m-k}\theta_1\cdots\theta_{\m},
\end{eqnarray}
we get that up to a scalar multiple,  ${\cal A}'_k$ has only one
highest weight vector $x_n^{\m-k}\theta_1\cdots\theta_{\m}$ and thus
it is irreducible.

It can be similarly proved that ${\cal A}'_k=\la
x_{2n}^{k-\m}\theta_1\cdots\theta_{\m}\ra$ is irreducible when
$k>\m$.

b) Assume $k=\m$. We claim that for any nonzero submodule $V$ of
${\cal A}'_{\m}$, we have
\begin{equation}
\theta_1\cdots\theta_{\m}\in V\;\;\mbox{or
}(x_{n-1}x_{2n}-x_nx_{2n-1})\theta_1\cdots\theta_{\m}\in V.
\end{equation}

In fact, there should be at least one weight vector $f\in V$ such
that $osp(2\m,2n)^+_{{0}}(f)=0$. Thus we can assume that $f$ is of
the form
\begin{equation}(x_{n-1}x_{2n}-x_nx_{2n-1})^l\theta_1\cdots\theta_{\m-1}
\theta_{\m}^{l_{\m}}\theta_{2\m}^{l_{2\m}}\end{equation} with
$l,l_{\m},l_{2\m}\in\{0,1\}$ such that $l_{\m}l_{2\m}=0$ or
\begin{equation}ax_n^{\m-i}\theta_1\cdots\theta_{\m}\theta_{\m+i+1}\cdots\theta_{2\m}+bx_{2n}^{\m-i}
\theta_1\cdots\theta_i\end{equation} with $0\leq i<\m$ and
$a,b\in\mathbb{C}$.

i) If
\begin{equation}f=ax_n^{\m-i}\theta_1\cdots\theta_{\m}\theta_{\m+i+1}\cdots\theta_{2\m}+bx_{2n}^{\m-i}\theta_1\cdots
\theta_i,\end{equation} then
\begin{eqnarray}& &
f_1=\prod\limits_{j=i+1}^{\m-1}(E_{j,2\m+2n}+E_{2\m+n,\m+j})(f)\nonumber\\
&=&a'x_n\theta_1\cdots\theta_{\m}\theta_{2\m}+b'x_{2n}\theta_1\cdots\theta_{\m-1}\in\la
f\ra.\end{eqnarray}
 If $a'\neq b'$, then
\begin{equation}
(E_{\m,2\m+2n}+E_{2\m+n,2\m})(f_1)=(-1)^{\m}(a'-b')\theta_1\cdots\theta_{\m}\in
V.\end{equation} Otherwise, $a'=b'$ and
\begin{equation}
(E_{\m,2\m+2n-1}+E_{2\m+n-1,2\m})(f_1)=(-1)^{\m}a'(x_{n-1}x_{2n}-x_nx_{2n-1})\theta_1\cdots\theta_{\m}\in
V.
\end{equation}

ii) If $f=\theta_1\cdots\theta_{\m-1}\theta_{2\m}$, then
\begin{eqnarray}
&&(E_{\m,2\m+2n-1}+E_{2\m+n-1,2\m})(E_{\m,2\m+n}-E_{2\m+2n,2\m})(f)\nonumber\\
&=&(-1)^{\m}(E_{\m,2\m+2n-1}+E_{2\m+n-1,2\m})(x_n\theta_1\cdots\theta_{\m}\theta_{2\m}+x_{2n}\theta_1\cdots\theta_{\m-1})\nonumber\\
&=&(x_{n-1}x_{2n}-x_nx_{2n-1})\theta_1\cdots\theta_{\m}\in V.
\end{eqnarray}

iii) When
$f=(x_{n-1}x_{2n}-x_nx_{2n-1})\theta_1\cdots\theta_{\m-1}\theta_{2\m}$,
we have
\begin{eqnarray}
&&(E_{\m,2\m+2n}+E_{2\m+n,2\m})(E_{\m,2\m+2n-1}+E_{2\m+n-1,2\m})(f)\nonumber\\
&=&(-1)^{\m}(E_{\m,2\m+2n}+E_{2\m+n,2\m})(x_n\theta_1\cdots\theta_{\m}\theta_{2\m}-x_{2n}\theta_1\cdots\theta_{\m-1})\nonumber\\
&=&-2\theta_1\cdots\theta_{\m}\in V.
\end{eqnarray}
Using Lemma 3.1, we get
\begin{equation}
x^\alpha\theta_1\cdots\theta_{\m}\in\la
\theta_1\cdots\theta_{\m}\ra+\la(x_{n-1}x_{2n}-x_nx_{2n-1})\theta_1\cdots\theta_{\m}\ra
\end{equation}
for $\alpha\in\mathbb{N}^{2n}$ such that
$\sum\limits_{i=1}^{\m}\alpha_i=\sum\limits_{i=1}^{\m}\alpha_{\m+i}$.
Now it is straightforward to check
\begin{equation}
{\cal A}'_{\m}\subset\la
\theta_1\cdots\theta_{\m}\ra+\la(x_{n-1}x_{2n}-x_nx_{2n-1})\theta_1\cdots\theta_{\m}\ra.
\end{equation}
\end{proof}

We can get basis for $\la\theta_1\cdots\theta_{\m}\ra$ and $\la
(x_{n-1}x_{2n}-x_nx_{2n-1})\theta_1\cdots\theta_{\m}\ra$ by the
following way. Set
\begin{equation}
\Delta=-\sum\limits_{i=1}^nx_i\partial_{x_{n+i}}+\sum\limits_{j=1}^{\m}\partial_{\theta_j}\partial_{\theta_{\m+j}}
\end{equation}
and
\begin{equation}
\eta=\sum\limits_{i=1}^nx_{n+i}\partial_{x_i}+\sum\limits_{j=1}^{\m}\theta_j\theta_{\m+j}.
\end{equation}
Let $L$ be the subalgebra of $osp(2\m,2n)$ consisting of the
matrices of the form
\begin{equation}
\left(
\begin{array}{cccc}
A & 0 & 0 & H_1\\
0 & -A^T & J & 0\\
0 &H_1^T & D & 0\\
-J^T & 0 & 0 & -D^T
\end{array}
\right).
\end{equation}
It is straightforward to check
\begin{equation}
g\Delta=\Delta g,\ \ g\eta=\eta g\ \ \mbox{for }g\in L.
\end{equation}
Set
\begin{eqnarray}
&&{\cal
A}_{s,t}=\mbox{Span}\;\{x^\alpha\prod\limits_{j=1}^{2\m}\theta_j^{l_j}\mid
\alpha\in\mathbb{N}^{2n};\;l_j\in\{0,1\};\;\sum\limits_{j=1}^{\m}l_{\m+j}-\sum\limits_{i=1}^n\alpha_i=s,\nonumber\\
&&\sum\limits_{j=1}^{\m}l_j+\sum\limits_{i=1}^n\alpha_{n+i}=t\}
\end{eqnarray}
for $s\in\mathbb{Z}$, $t\in\mathbb{N}$ and
\begin{equation}
{\cal H}_{s,t}=\{f\in{\cal A}_{s,t}\mid \Delta(f)=0\}.
\end{equation}
We have
\begin{equation}
{\cal A}_{s,t}={\cal H}_{s,t}\oplus\eta{\cal
A}_{s-1,t-1}=\bigoplus_{l=0}^t\eta^{t-l}{\cal H}_{s-t+l,l}.
\end{equation}
Thus
\begin{equation}
{\cal A}_{\m}'=\bigoplus\limits_{l=0}^\infty{\cal
A}_{\m-l,l}=\bigoplus\limits_{l=0}^\infty\bigoplus\limits_{q=0}^{l}\eta^{l-q}{\cal
H}_{\m-2l+q,q}=\bigoplus_{t,q=0}^{\infty}\eta^{t}{\cal
H}_{\m-2t-q,q}.
\end{equation}
By  similar arguments as those in 1) of the proof of Theorem 2.2, we
get
\begin{equation}
{\cal H}_{\m-2t-q,q}=\left\{
\begin{array}{lll}
U(L)\big(x_n^{2t}\theta_1\cdots\theta_q\theta_{\m+q+1}\cdots\theta_{2\m}\big)
&\mbox{if} &q\leq \m,\\
U(L)\big(x_n^{2t}(x_{n-1}x_{2n}-x_nx_{2n-1})^{q-\m}\theta_1\cdots\theta_{\m}\big)
& \mbox{if} & q>\m.
\end{array}
\right.
\end{equation}
Since
\begin{eqnarray}
&&\eta^t(x_n^{2t}\theta_1\cdots\theta_q\theta_{\m+q+1}\cdots\theta_{2\m})=\frac{(2t)!}{t!}x_n^tx_{2n}^t\theta_1\cdots\theta_q\theta_{\m+q+1}\cdots\theta_{2\m}\nonumber\\
&=&\frac{(2t)!}{t!}(-E_{2n,n})^t(\theta_1\cdots\theta_q\theta_{\m+q+1}\cdots\theta_{2\m})\nonumber\\
&\in& \left\{
\begin{array}{lll}
\la \theta_1\cdots\theta_{\m}\ra &\mbox{if} & \m-q\mbox{ is even},\\
\la \theta_1\cdots\theta_{\m-1}\theta_{2\m}\ra=\la
(x_{n-1}x_{2n}-x_nx_{2n-1})\theta_1\cdots\theta_{\m}\ra &\mbox{if} &
\m-q\mbox{ is odd}
\end{array}
\right.
\end{eqnarray}
when $q\leq \m$ and
\begin{eqnarray}
&&\eta^t(x_n^{2t}(x_{n-1}x_{2n}-x_nx_{2n-1})^{q-\m}\theta_1\cdots\theta_{\m})\nonumber\\
&=&
\frac{(2t)!}{t!}x_n^tx_{2n}^t(x_{n-1}x_{2n}-x_nx_{2n-1})^{q-\m}\theta_1\cdots\theta_{\m}\nonumber\\
&\in& \left\{
\begin{array}{lll}
U(sp(2n,\mathbb{C}))(1)\theta_1\cdots\theta_{\m} &\mbox{if} & \m-q\mbox{ is even},\\
U(sp(2n,\mathbb{C}))(x_{n-1}x_{2n}-x_nx_{2n-1})\theta_1\cdots\theta_{\m}
&\mbox{if} & \m-q\mbox{ is odd}
\end{array}
\right.
\end{eqnarray}
when $q>\m$ (cf. [\ref{r6}]), we get
\begin{equation}
\eta^t{\cal H}_{\m-2t-q,q}\subset \left\{
\begin{array}{lll}
\la \theta_1\cdots\theta_{\m}\ra &\mbox{if} & \m-q\mbox{ is even},\\
\la (x_{n-1}x_{2n}-x_nx_{2n-1})\theta_1\cdots\theta_{\m}\ra
&\mbox{if} & \m-q\mbox{ is odd}.
\end{array}
\right.
\end{equation}
Therefore,
\begin{equation}
\la\theta_1\cdots\theta_{\m}\ra=\bigoplus_{q\in\mathbb{N};\;\m-q\;\rm{is\;
even}}\eta^t{\cal H}_{\m-2t-q,q}
\end{equation}
and
\begin{equation}
\la
(x_{n-1}x_{2n}-x_nx_{2n-1})\theta_1\cdots\theta_{\m}\ra=\bigoplus_{q\in\mathbb{N};\;\m-q\;\rm{is\;
odd}}\eta^t{\cal H}_{\m-2t-q,q}.
\end{equation}
Like in Theorem 2.3, we denote
\begin{eqnarray}
h(\vec{k},\vec{l},\vec{s})&=&\prod\limits_{t=1}^{n}
x_t^{k_t}\prod\limits_{1\leq i<j\leq
n}(x_ix_{n+j}-x_jx_{n+i})^{k_{i,j}}\prod\limits_{j=1}^{\m}
\theta_j^{l_j}\theta_{\m+j}^{l_{\m+j}}\prod\limits_{1\leq i<j\leq
\m}(\theta_i\theta_{\m+i}\nonumber\\
&&-\theta_j\theta_{\m+j})^{l_{i,j}}\prod\limits_{1\leq p\leq n,1\leq
q\leq \m}(x_{n+p}-x_p\theta_q\theta_{\m+q})^{s_{p,q}}
\end{eqnarray}
where
\begin{eqnarray}
&&\vec{k}=(k_1,\cdots,k_{n};k_{1,2},k_{1,3},\cdots,k_{1,n},k_{2,3}\cdots,k_{n-1,n})\in\mathbb{N}^{\frac{n(n+1)}{2}},\\
&&\vec{l}=(l_1,\cdots,l_{\m};l_{1,2},\cdots,l_{1,\m},l_{2,3},\cdots,l_{2,\m},\cdots,l_{\m-1,\m})\in\{0,1\}^{\frac{\m(\m+1)}{2}},\\
&&\vec{s}=(s_{1,1},\cdots,s_{1,\m},\cdots,s_{n,\m})\in\{0,1\}^{\m
n}.
\end{eqnarray}
Set
\begin{eqnarray}
I&=&\{(\vec{k},\vec{l},\vec{s})\mid l_t+l_{\m+t}+\sum\limits_{1\leq
i<t}l_{i,t}+\sum\limits_{t<j\leq
\m}l_{t,j}+\sum\limits_{p=1}^{n}s_{p,t}\leq 1\mbox{ for
}t\in\overline{1,\m};\nonumber\\
&&k_{i,j}k_t=0\mbox{ for }i<j<t; \; k_{i,j}k_{i',j'}=0\mbox{ for
}i>i'\mbox{ and }j<j';\nonumber\\
&&k_tl_{i,j}=0\mbox{ for }t\in\overline{1,n},1\leq i<j\leq \m;
\;k_ts_{p,q}=0\mbox{ for }t<p;\;\nonumber\\
&& k_{i,j}s_{p,q}=0\mbox{ for }i<j<p;\;s_{p,q}s_{p',q'}=0\mbox{ for
}p>p'\mbox{ and }q<q';\nonumber\\
&& l_{i,j}=0 \mbox{ if
}l_t=l_{\m+t}=\sum\limits_{p=1}^{n}s_{p,t}=\sum\limits_{i'<t}l_{i',t}+\sum\limits_{j'>t}l_{t,j'}=0\mbox{
for some }i<t<j; \nonumber\\
&&l_{i,j}l_{i',j'}=0\mbox{ if }i<i'<j<j';\; l_{i,j}s_{p,q}=0\mbox{
if }i<j<q \}.
\end{eqnarray}
Then the subspace ${\cal H}_{\m-2t-q,q}$ has a basis
\begin{eqnarray}
&&B_{t,q}=\{h(\vec{k},\vec{l},\vec{s})\mid
\;\sum_{j=1}^{\m}l_{\m+j}+\sum_{1\leq
i<j\leq\m}l_{i,j}-\sum_{i=1}^nk_i-\sum_{1\leq i\leq j\leq
n}k_{i,j}=\m-2t-q;\nonumber\\
&&\;\sum_{j=1}^{\m}l_{j}+\sum_{1\leq
i<j\leq\m}l_{i,j}+\sum_{p=1}^n\sum_{j=1}^{\m}s_{p,j}+\sum_{1\leq
i\leq j\leq n}k_{i,j}=q;\;(\vec{k},\vec{l},\vec{s})\in I\}.
\end{eqnarray}
Note
\begin{equation}
\{f\in{\cal A}_{s,t}\mid \eta(f)=0\}=0
\end{equation}
when $s+t<\m$. Hence we have:

\begin{thm} The set
\begin{equation}
\bigcup_{t,q\in\mathbb{N};\m-q\rm{\; is\; even}}\eta^t B_{t,q}
\end{equation}
forms a basis for $\la\theta_1\cdots\theta_{\m}\ra$ and the set
\begin{equation}
\bigcup_{t,q\in\mathbb{N};\m-q\rm{\; is\; odd}}\eta^t B_{t,q}
\end{equation}
forms a basis for
$\la(x_{n-1}x_{2n}-x_nx_{2n-1})\theta_1\cdots\theta_{\m}\ra$.
\end{thm}

\section{Proof of Theorem 3}
\setcounter{equation}{0}

 Recall the Lie superalgebra
\begin{equation}
osp(2\m+1,2n)=osp(2\m+1,2n)_{0}\oplus osp(2\m+1,2n)_{1}
\end{equation}
where
\begin{eqnarray}
&&osp(2\m+1,2n)_{0}=\sum\limits_{i,j=1}^{\m}\big(\mathbb{C}(E_{i,j}-E_{\m+j,\m+i})+\mathbb{C}(E_{i,\m+j}-E_{j,\m+i})\nonumber\\
&&+\mathbb{C}(E_{\m+i,j}-E_{\m+j,i})\big)+\sum_{i=1}^{\m}\big(\mathbb{C}(E_{i,2\m+1}-E_{2\m+1,\m+i})+\mathbb{C}(E_{\m+i,2\m+1}-E_{2\m+1,i})\big)\nonumber\\
&&+\sum\limits_{p,q=1}^n\big(\mathbb{C}(E_{2\m+1+p,2\m+1+q}-E_{2\m+1+n+q,2\m+1+n+p})+\mathbb{C}(E_{2\m+1+p,2\m+1+n+q}\nonumber\\
&&+E_{2\m+1+q,2\m+1+n+p})+\mathbb{C}(E_{2\m+1+n+p,2\m+1+q}+E_{2\m+1+n+q,2\m+1+p})\big)
\end{eqnarray}
and
\begin{eqnarray}
&&osp(2\m+1,2n)_{1}=\sum\limits_{i\in\overline{1,\m};p\in\overline{1,n}}\big(\mathbb{C}(E_{i,2\m+1+p}
-E_{2\m+1+n+p,\m+i})\nonumber\\&&+\mathbb{C}(E_{i,2\m+1+n+p}
+E_{2\m+1+p,\m+i})+\mathbb{C}(E_{\m+i,2\m+1+p}-E_{2\m+1+n+p,i})\nonumber\\
&&+\mathbb{C}(E_{\m+i,2\m+1+n+p}+E_{2\m+1+p,\m+i})\big)
+\sum_{p=1}^n\big(\mathbb{C}(E_{2\m+1,2\m+1+p}-E_{2\m+1+n+p,2\m+1})\nonumber\\
&&+\mathbb{C}(E_{2\m+1,2\m+1+n+p}+E_{2\m+1+p,2\m+1})\big).
\end{eqnarray}
Take
\begin{equation}
H=\sum\limits_{i=1}^{\m}\mathbb{C}(E_{i,i}-E_{\m+i,\m+i})+\sum\limits_{j=1}^n\mathbb{C}(E_{2\m+1+j,2\m+1+j}-E_{2\m+1+n+j,2\m+1+n+j})
\end{equation}
as a Cartan subalgebra of $osp(2\m+1,2n)$. We still denote by
$\lambda_1,\cdots,\lambda_{\m}$, $\nu_1,\cdots,\nu_n$ the
fundamental weights. Let
\begin{eqnarray}
osp(2\m+1,2n)^+&=&\sum\limits_{1\leq i<j\leq
m}\big(\mathbb{C}(E_{i,j}-E_{\m+j,\m+i})+\mathbb{C}(E_{i,\m+j}-E_{j,\m+i})\big)\nonumber\\
&&+\sum\limits_{1\leq p<q\leq
n}\mathbb{C}(E_{2\m+1+p,2\m+1+q}-E_{2\m+1+n+q,2\m+1+n+p})\nonumber\\
&&+\sum\limits_{1\leq p\leq q\leq
n}\mathbb{C}(E_{2\m+1+p,2\m+1+n+q}+E_{2\m+1+q,2\m+1+n+p})\nonumber\\
&&+\sum\limits_{1\leq i\leq m,1\leq q\leq
n}\big(\mathbb{C}(E_{i,2\m+1+q}-E_{2\m+1+n+q,\m+i})\nonumber\\
&&+\mathbb{C}(E_{i,2\m+1+n+q}+E_{2\m+1+q,\m+i})\big)\nonumber\\
&&+\sum_{i=1}^{\m}\mathbb{C}(E_{i,2\m+1}-E_{2\m+1,\m+i})\nonumber\\&&+\sum_{p=1}^n\mathbb{C}(E_{2\m+1,2\m+1+n+p}+E_{2\m+1+p,2\m+1}),
\end{eqnarray}
and $osp(2\m+1,2n)^+_{\sigma}=osp(2\m+1,2n)_{\sigma}\cap
osp(2\m+1,2n)^+$ for $\sigma=0,1$. We redefine $\Delta$ and $\eta$
by
\begin{eqnarray}
&&\Delta_x=-2\sum\limits_{i=1}^r
x_i\partial_{x_{\m+i}}+2\sum\limits_{i=r+1}^{\m}\partial_{x_i}\partial_{x_{\m+i}}+\partial^2_{x_{2\m+1}},\\
&&\Delta_\theta=\sum\limits_{j=1}^n\partial_{\theta_j}\partial_{\theta_{n+j}},\;\;\Delta=\Delta_x+2\Delta_\theta,\\
&&\eta_x=
2\sum\limits_{i=1}^rx_{\m+i}\partial_{x_i}+2\sum\limits_{i=r+1}^{\m}x_ix_{\m+i}+x_{2\m+1}^2,\\
&&\eta_\theta=\sum\limits_{j=1}^n\theta_j\theta_{n+j},\;\;\eta=\eta_x+2\eta_\theta.
\end{eqnarray}
Recall the notions in (1.7) and (1.10).

\begin{thm}
The subspace ${\cal H}^r_k$ ($0\leq r\leq \m$) is an irreducible
highest weight $osp(2\m+1,2n)$-module with highest weights and
corresponding vectors are listed as follows:
\begin{center}
\begin{tabular}{c|c|c|c}
  \hline
  $r$ & $k$ & vector & weight\\
  \hline
  $r=0$ & $k>0$ & $x_1^k$ & $k\lambda_1$\\
 \hline
 $r<\m$ & $k>0$ & $x_{r+1}^k$ & $-(k+1)\lambda_r+k\lambda_{r+1}$\\
 \cline{2-4}
 & $k\leq0$ & $x_r^{-k}$ & $k\lambda_{r-1}-(k+1)\lambda_r$\\
 \hline
 $r=\m$ & $k\leq0$ & $x_{\m}^{-k}$ & $k\lambda_{\m-1}-(k-1)\lambda_{\m}$\\
 \cline{2-4}
  & $0<k\leq n$ & $\theta_1\cdots\theta_k$ & \begin{tabular}{lll}$\nu_k$ & if & $k<n$\\
   $2\nu_n$ & if & $k=n$
   \end{tabular}\\
 \cline{2-4}
  & $n<k\leq 2n$ &
  $\sum\limits_{i=0}^{k-n-1}\frac{(k-n-1-i)!}{2^ii!(2k-2n-1-2i)!}x_{2\m+1}^{2k-2n-1-2i}\eta_{\theta}^i\theta_1\cdots\theta_{k-n-1}$
  &\begin{tabular}{lll}
   $2\nu_n$ & if & $k=n+1$\\
   $\nu_{k-n-1}$ & if & $k>n+1$\\
   \end{tabular}\\
 \cline{2-4}
 & $k>2n$
 &$\begin{array}{l}\sum\limits_{j=0}^n\sum\limits_{i=0}^{k-2n-1}\frac{1}{i!(k-2n-1-i)!(2k-2n-1-2i-2j)!!j!}\\ \times x_{\m}^{k-2n-1-i}x_{2\m}^ix_{2\m+1}^{2k-2n-1-2i-2j}\eta_{\theta}^j\end{array}$
 & $\begin{array}{l}(2n-k+1)\lambda_{\m-1}\\ +(k-2n-2)\lambda_{\m}\end{array}$\\
 \hline
\end{tabular}
\end{center}
Moreover, ${\cal A}^r_k={\cal H}^r_k\oplus\eta{\cal A}^r_{k-2}$.
\end{thm}
\begin{proof}
(1) Assume $r=0$.

First we can get $H^0_k=\la x_1^k\ra$ by induction on $n$. If
\begin{equation}
osp(2\m+1,2n)_{0}^+(g)=0,\label{b20}
\end{equation}
then up to a scalar multiple, $g$ must be of the  form
\begin{equation}\sum\limits_{i=0}^l
a_ix_1^{k-2l-t}\eta_x^{l-i}\eta_\theta^i\theta_1\cdots\theta_t\end{equation}
with $l-2l-t\geq0$, $l\geq0$, $0\leq t\leq n$ and $a_i=0$ for
$i>n-t$. Since
\begin{equation}
\begin{array}{lll}
0&=&\Delta(\sum\limits_{i=0}^l
a_ix_1^{k-2l-t}\eta_x^{l-i}\eta^i\theta_1\cdots\theta_t)\\
&=&2\sum\limits_{i=0}^{l-1}
a_i(l-i)(2k-2l-2t+2\m-1-2i)x_1^{k-2l-t}\eta_x^{l-i-1}\eta^i\theta_1\cdots\theta_t\\
&&-2\sum\limits_{i=1}^l
a_ii(n-t-i+1)x_1^{k-2l-t}\eta_x^{l-i}\eta^{i-1}\theta_1\cdots\theta_t\\
\end{array}\end{equation}
we get that
\begin{equation}a_i(l-i)(2k-2l-2t+2\m-1-2i)=a_{i+1}(i+1)(n-t-i)\end{equation}
for $0\leq i<l$  and $l\leq n-t$. Thus all the vectors satisfying
(\ref{b20}) should be a scalar multiple of
\begin{equation}
f_{l,t}=\sum\limits_{i=0}^l
\frac{(n-t-i)!}{i!(l-i)!(2k-2l-2t+2\m-1-2i)!}x_1^{k-2l-t}\eta_x^{l-i}\eta_\theta^i\theta_1\cdots\theta_t,\label{d5}
\end{equation}
where $0\leq t\leq n$ and $0\leq
l\leq\mbox{min}\{n-t,\frac{1}{2}(k-t)\}$.

Take any $0\neq f\in{\cal H}^{0}_k$. Then there should be some
$f_{l_0,t_0}\in\la f\ra$. If $l_0=0$, then
\begin{equation}
x_1^k=\prod\limits_{i=1}^t(-1)^{t-1}(E_{1,2\m+1+i}-E_{2\m+1+n+i,\m+1})(f_{0,t_0}),
\end{equation}
and so ${\cal H}^{0}_k=\la x_1^k\ra\subset\la f_{0,t_0}\ra\subset\la
f\ra$. Now assume $l_0>0$. Observe that
\begin{equation}
(E_{1,2\m+n+t+1}+E_{2\m+1+t+1,\m+1})(f_{l,t})=(-1)^{t-1}\big(2(n-t-k+l+t-m)+1\big)f_{l-1,t+1}\label{b6}
\end{equation}
for $0<l\leq\mbox{min}\{n-t,\frac{1}{2}(k-t)\}$. Since
$2(n-t-k+l+t-m)+1$ is odd, we obtain
$f_{l_0-1,t_0+1},\cdots,f_{0,l_0+t_0}\in\la f\ra$ by (\ref{b6}).
Hence ${\cal H}^{0}_k=\la f_{0,l_0+t_0}\ra\subset\la f\ra$.\psp

(2) Suppose $r<\m$.

Taking induction on $n$, we obtain that the submodule ${\cal H}^r_k$
is generated by $x_{r+1}^k$ if $k\geq 0$ or $x_r^{-k}$ if $k<0$.

Let $W$ be a nonzero submodule of ${\cal H}^r_k$. Then $W$ should
contain some weight vector $f$ annihilated by $osp(2\m+1,2n)^+_0$.
Indeed, $f$ should be of the form
\begin{equation}
\sum_{i=0}^la_i\eta_x^{l-i}x_{r+1}^{k-2l-t}\eta_{\theta}^i\theta_1\cdots\theta_t
\end{equation}
or
\begin{equation}
\sum_{i=0}^la_i\eta_x^{l-i}x_r^{2l+t-k}\eta_{\theta}^i\theta_1\cdots\theta_t
\end{equation}
because $osp(2\m+1,2n)_0=so(2\m+1)\oplus sp(2n)$, where $0\leq t\leq
n$ and $0\leq l\leq n-t$.

If
\begin{equation}
f=\sum_{i=0}^la_i\eta_x^{l-i}x_{r+1}^{k-2l-t}\eta_{\theta}^i\theta_1\cdots\theta_t,
\end{equation}
then by the similar argument as in (1), we obtain $x_{r+1}^k\in\la
f\ra\subset W$, which implies $W={\cal H}^r_k$. Now we assume
\begin{equation}
f=\sum_{i=0}^la_i\eta_x^{l-i}x_r^{2l+t-k}\eta_{\theta}^i\theta_1\cdots\theta_t.
\end{equation}
Since
\begin{eqnarray}
0=\Delta(f)&=&\sum_{i=0}^{l-1}2(l-i)(2l-2p+2\m-2r-1-2i)a_i\eta_x^{l-i-1}x_r^p\eta_{\theta}^i\theta_1\cdots\theta_t\nonumber\\
&&-\sum_{i=1}^l2i(n-t-i+1)a_i\eta_x^{l-i}x_r^p\eta_{\theta}^{i-1}\theta_1\cdots\theta_t,
\end{eqnarray}
we get
\begin{equation}
a_{i+1}(i+1)(n-t-i)=a_i(l-i)(2l-2p+2\m-2r-1-2i).
\end{equation}
So $f$ is a scalar multiple of
\begin{equation}
f_{l,t}=\sum_{i=0}^l\frac{(n-t-i)!(2l-2p+2\m-2r-1)!!}{i!(l-i)!(2l-2p+2\m-2r-1-2i)!!}\eta_x^{l-i}x_r^{2l+t-k}\eta_{\theta}^i\theta_1\cdots\theta_t,
\end{equation}
where
\begin{equation}
\frac{(2l-2p+2\m-2r-1)!!}{(2l-2p+2\m-2r-1-2i)!!}=\prod_{j=0}^{i-1}(2l-2p+2\m-2r-1-2j).
\end{equation}
Note
\begin{eqnarray}
&&(E_{r,2\m+1+n+t+1}+E_{2\m+1+t+1,\m+r})(f_{l,t})\nonumber\\
=&&(-1)^{t+1}(2l+t-k)(2l-2p+2m-2r-1-2n+2t)f_{l-1,t+1}
\end{eqnarray}
for $2l+t>k$, $t<n$ and $l>0$. Thus we have
$f_{l-1,t+1},f_{l-2,t+2},\cdots,f_{0,l+t}\in W$ when $k\leq l+t\leq
n$ or $f_{l-1,t+1},f_{l-2,t+2},\cdots,f_{k-l-t,2(l+t)-k}\in W$ when
$l+t<k$. For the later case, we can get $x_{r+1}^k\in W$. Now
suppose $f_{0,l+t}\in W$. Applying
\begin{equation}
E_{r,2\m+1+q}-E_{2\m+1+n+q,\m+r})|_{{\cal
A}^r_k}=\partial_{x_r}\partial_{\theta_q}-\theta_{n+q}\partial_{x_{\m+r}}
\end{equation}
to $f_{0,l+t}$, we obtain $x_r^{-k}\in W$ or
$\theta_1\cdots\theta_k\in W$. Applying
\begin{equation}
(E_{r+1,2\m+1+q}-E_{2\m+1+n+q,\m+r+1})|_{{\cal
A}^r_k}=x_{r+1}\partial_{\theta_q}-\theta_{n+q}\partial_{x_{\m+r+1}}
\end{equation}
to $\theta_1\cdots\theta_k$, we get $x_{r+1}^k\in W$. Thus $W={\cal
H}^r_k$.

(3) Assume $r=\m$.

Denote
\begin{equation}
g_{p,q}=\sum_{s=0}^q\frac{1}{s!(q-t)!(p+q-2s)!!}x_{2\m+1}^{p+q-2s}x_\m^{q-t}x_{2\m}^s
\end{equation}
for $p>q\geq0$. By induction on $n$, we can obtain
\begin{equation}
{\cal H}^{\m}_k=\left\{
\begin{array}{ll}
\la x_\m^{n-k}\theta_1\cdots\theta_n\ra &\mbox{if }k\leq n,\\
\la g_{k-n,k-n-1}\theta_1\cdots\theta_n\ra &\mbox{if
}k>n.\label{b114}
\end{array}
\right.
\end{equation}
Let $W$ be a nonzero submodule of ${\cal H}^r_k$. Then $W$ should
contain some weight vector $f$ annihilated by $osp(2\m+1,2n)^+_0$.
Indeed, $f$ should be of the form
\begin{equation}
\sum_{i=0}^l
a_i\eta_x^{l-i}x_\m^{2l+t-k}\eta_\theta^i\theta_1\cdots\theta_t
\end{equation}
or
\begin{equation}
\sum_{i=0}^{l'}b_i\eta_x^{l'-i}g_{k-2l'-t',k-2l'-t'-1}\eta_\theta^i\theta_1\cdots\theta_{t'}
\end{equation}
because $osp(2\m+1,2n)_0=so(2m_1+1)\oplus sp(2n)$, where $0\leq
t,t'\leq n$, $0\leq l\leq n-t$, $0\leq l'\leq n-t'$, $2l+t\geq k$
and $2l'+t'<k$.

(i) Suppose
$f=\sum_{i=0}^{l'}b_i\eta_x^{l'-i}g_{k-2l'-t',k-2l'-t'-1}\eta_\theta^i\theta_1\cdots\theta_{t'}$.
Since
\begin{eqnarray}
0=\Delta(f)&=&\sum_{i=0}^{l'-1}2b_i(2k-2t'-2l'-1-2i)(l-i)\eta_x^{l'-i-1}g_{k-2l'-t',k-2l'-t'-1}\eta_\theta^i\theta_1\cdots\theta_{t'}\nonumber\\
&&-\sum_{i=1}^{l'}2b_ii(n-t'-i+1)\eta_x^{l'-i}g_{k-2l'-t',k-2l'-t'-1}\eta_\theta^{i-1}\theta_1\cdots\theta_{t'},
\end{eqnarray}
we get
\begin{equation}
b_{i+1}(i+1)(n-t'-i)=b_i(l'-i)(2k-2t'-2l'-1-2i).
\end{equation}
Thus $f$ should be a scalar multiple of
\begin{equation}
h_{l,t}=\sum_{i=0}^l\frac{(n-t-i)!(2k-2t-2l-1)!!}{i!(l-i)!(2k-2t-2l-1-2i)!!}
\eta_x^{l-i}g_{k-2l-t,k-2l-t-1}\eta_\theta^i\theta_1\cdots\theta_t.
\end{equation}
Note
\begin{equation}
(E_{\m,2\m+1+n+t+1}+E_{2\m+1+t+1,2\m})(h_{l,t})=\frac{(-1)^t(n-t-l)}{(2k-4l-2t-1)(2k-2t-2l+1)}h_{l,t+1}
\end{equation}
if $l+t<n$ and $k-2l-t>1$. So we can assume $l+t=n$ or $k-2l-t=1$.

(a) If $l+t=n$, then
\begin{equation}
h_{n-t,t}=\sum_{i=0}^{n-t}\frac{(2k-2n-1)!!}{i!(2k-2n-1-2i)!!}
\eta_x^{n-t-i}g_{k-2n+t,k-2n+t-1}\eta_\theta^i\theta_1\cdots\theta_t.
\end{equation}
We have
\begin{eqnarray}
&
&(E_{2\m,2\m+1+n+t+1}+E_{2\m+1+t+1,\m})(h_{n-t,t})\nonumber\\&=&(-1)^{t+1}(k-2n+t)(2k-4n+2t+1)h_{n-t-1,t+1}.
\end{eqnarray}
Thus $h_{n-t-1,t+1},h_{n-t-2,t+2},\cdots,h_{0,n}\in W$. Since
$h_{0,n}=g_{k-n,k-n-1}\theta_1\cdots\theta_n$, we get $W={\cal
H}^{\m}_k$ by (\ref{b114}).

(b) Assume $k-2l-t=1$ and $l+t<n$. Then
\begin{equation}
h_{l,t}=h_{l,k-2l-1}=\sum_{i=0}^l\frac{(n-k+2l+1-i)!(2l-1)!!}{i!(l-i)!(2l-1-2i)!!}
x_{2\m+1}^{2l-2i+1}\eta_\theta^i\theta_1\cdots\theta_{k-2l-1}.
\end{equation}
Now we take induction on $l$. When $l=0$, we have
\begin{equation}
(E_{2\m+1,2\m+1+n+k}+E_{2\m+1+k,2\m+1})(h_{0,k-1})=(-1)^{k-1}\theta_1\cdots\theta_k
\end{equation}
and
\begin{equation}
x_\m^{n-k}\theta_1\cdots\theta_n=(-1)^{(k+1)(n-k)}\prod_{j=k+1}^n(E_{2m,2\m+1+n+j}+E_{2\m+1+j,\m})(\theta_1\cdots\theta_k).
\end{equation}
Thus $W={\cal H}^{\m}_k$.

Now we assume $l>0$. Observe that $k-2l=t+1\leq t+l<n$ and
\begin{eqnarray}
&&(E_{2\m+1,2\m+1+n+k-2l+1}+E_{2\m+1+k-2l+1,2\m+1})(E_{2\m+1,2\m+1+n+k-2l}\nonumber\\
&&+E_{2\m+1+k-2l,2\m+1})(h_{l,k-2l-1})\nonumber\\
&=&-(n-k+l+1)(2l-1)(2n-2k+2l-1)h_{l-1,k-2l+1}.
\end{eqnarray}
Since $n-k+l+1=n-l-t>0$, we get $h_{l-1,k-2l+1}\in W$. Therefore
$W={\cal H}^{\m}_k$ by inductive assumption.

(ii) Assume $f=\sum_{i=0}^l
a_i\eta_x^{l-i}x_\m^{2l+t-k}\eta_\theta^i\theta_1\cdots\theta_t$.
Then
\begin{eqnarray}
0=\Delta(f)&=&\sum_{i=0}^{l-1}2a_i(2k-2t-2l-1-2i)(l-i)\eta_x^{l-i-1}x_{\m}^{2l+t-k}\eta_\theta^i\theta_1\cdots\theta_t\nonumber\\
&&-\sum_{i=1}^l2a_ii(n-t-i+1)\eta_x^{l-i}x_\m^{2l+t-k}\eta_\theta^{i-1}\theta_1\cdots\theta_t.
\end{eqnarray}
Thus we get that $f$ should be a scalar multiple of
\begin{equation}
f_{l,t}=\sum_{i=0}^l\frac{(n-t-i)!(2k-2t-2l-1)!!}{i!(l-i)!(2k-2t-2l-1-2i)!!}
\eta_x^{l-i}x_\m^{2l+t-k}\eta_\theta^i\theta_1\cdots\theta_t,
\end{equation}
where we treat
\begin{equation}
\frac{(2k-2t-2l-1)!!}{(2k-2t-2l-1-2i)!!}=(2k-2t-2l-1)(2k-2t-2l-3)\cdots(2k-2t-2l+1-2i).
\end{equation}

We will show $W={\cal H}^{\m}_k$ by induction on $l$. When $l=0$, we
have $f_{0,t}=x_{\m}^{t-k}\theta_1\cdots\theta_t\in W$ and
\begin{equation}
x_\m^{n-k}\theta_1\cdots\theta_n=(-1)^{(k+1)(n-k)}\prod_{j=k+1}^n(E_{2m,2\m+1+n+j}+E_{2\m+1+j,\m})(f_{l,k-2l}).
\end{equation}
Thus $W={\cal H}^{\m}_k$ because of (\ref{b114}). Now assume $l>0$.
Note
\begin{equation}
(E_{\m,2\m+1+n+t+1}+E_{2\m+1+t+1,2\m})(f_{l,t})=(-1)^t(2l+t-k)(2n-2k+2l+1)f_{l-1,t+1}.
\end{equation}
If $2l+t\neq k$, then we get $f_{l-1,t+1}\in W$, which implies
$W={\cal H}^{\m}_k$ by inductive assumption.

If $2l+t=k$, then
\begin{equation}
f_{l,t}=f_{l,k-2l}=\sum_{i=0}^l\frac{(2l-1)!!(n-k+2l-i)!}{i!(l-i)!(2l-1-2i)!!}x_{2\m+1}^{2l-2i}\eta_{\theta}^i\theta_1\cdots\theta_{k-2l}.
\end{equation}
When $k-2l=n-1$, we have $l=1$, $k=n+1$ and
\begin{eqnarray}
&&(E_{2\m+1,2\m+1+2n}+E_{2\m+1+n,2\m+1})(f_{1,n-1})\nonumber\\
=&&(-1)^kx_{2\m+1}\theta_1\cdots\theta_{n}.
\end{eqnarray}
Therefore, $W={\cal H}^{\m}_k$ due to (\ref{b114}).

Now we assume $k-2l<n-1$. Observe that
\begin{eqnarray}
&&(E_{2\m+1,2\m+1+n+k-2l+2}+E_{2\m+1+k-2l+2,2\m+1})(E_{2\m+1,2\m+1+n+k-2l+1}\nonumber\\
&&+E_{2\m+1+k-2l+1,2\m+1})(f_{l,k-2l})\nonumber\\
&=&-(2n-2k+2l+1)(2l-1)(n-k+l)f_{l-1,k-2l+2}.
\end{eqnarray}
We get $W={\cal H}^{\m}_k$ by inductive assumption if $n-k+l\neq0$.

When $k=n+l$, we have
\begin{equation}
f_{l,t}=f_{k-n,2n-k}=\sum_{i=0}^{k-n}\frac{(2k-2n-1)!!}{i!(2k-2n-1-2i)!!}x_{2\m+1}^{2(k-n-i)}\eta_\theta^i\theta_1\cdots\theta_{2n-k}
\end{equation}
and
\begin{equation}
(E_{2\m,2\m+1+3n-k+1}+E_{2\m+1+2n-k+1,\m})(f_{k-n,2n-k})=h_{k-n-1,2n-k+1}.
\end{equation}
Thus we get $h_{k-n-1,2n-k-1}\in W$, which implies $W={\cal
H}^{\m}_k$ by (i).
\end{proof}

Let us go for the submodule ${\cal A}'_k$ (cf. (1.15)). We know that
the submodule ${\cal A}'_k$ is irreducible when $n=0$. Thus it is
not difficult to check the irreducibility of ${\cal A}'_k$.

\begin{thm}The $osp(2\m+1,2n)$-submodules ${\cal A}'_k$ are
irreducible.
\end{thm}
\begin{proof}
(1) Assume $S_1\neq\emptyset$. Take $i\in S_1$ and $0\neq f\in{\cal
A}'_k$. If $f$ is not independent of $\theta_{2\m+1}$, we apply
\begin{equation}
E_{2\m+1+i,2\m+1+n+i}|_{{\cal A}'}=x_i\partial_{x_{n+i}}
\end{equation}
to $f$ and get some $0\neq f_1\in\la f\ra$ satisfying
$\partial_{x_{n+i}}(f_1)=0$. If $f_1$ is not independent of
$\theta_{2\m+1}$, we have
\begin{equation}
f_2=(E_{2\m+1,2\m+1+n+i}+E_{2\m+1+i,2\m+1})|_{{\cal
A}'}(f_1)=x_i\partial_{2\m+1}(f_1).
\end{equation}
Anyway there is some nonzero
$f'=f'(x_1,\cdots,x_{2n};\theta_1\cdots,\theta_{2\m})\in\la f\ra$.
By Theorem 3.2, we obtain
\begin{equation}
\mbox{Span}\;\{x^\alpha\theta_{i_1}\cdots\theta_{i_t}\in{\cal
A}'_k\mid 0\leq t\leq
2\m;\;i_1,\cdots,i_t\in\overline{1,2\m}\}\subset\la f\ra.
\end{equation}
Now for any
$x^\alpha\theta_{i_1}\cdots\theta_{i_t}\theta_{2\m+1}\in{\cal
A}'_k$, we have
\begin{equation}
x^\alpha\theta_{i_1}\cdots\theta_{i_t}\theta_{2\m+1}=\frac{(-1)^t}{\alpha_i+1}(E_{2\m+1,2\m+1+i}-E_{2\m+1+n+i,2\m+1})(x_ix^\alpha\theta_{i_1}\cdots\theta_{i_t})\in\la
f\ra.
\end{equation}
Therefore $\la f\ra={\cal A}'_k$.

It can be similarly proved when $T_1\neq\emptyset$.

(2) Suppose $S_1\cup T_1=\emptyset$. We can assume
$T=\overline{1,n}$ by symmetry. We claim that ${\cal A}'_k=\la
x_n^{m_1-k}\theta_1\cdots\theta_{m_1}\ra$ when $k\leq\m$.

Using Theorem 3.2, we get
\begin{equation}
\mbox{Span}\;\{x^\alpha\theta_{i_1}\cdots\theta_{i_t}\in{\cal
A}'_k\mid 0\leq
t\leq2\m;\;i_1,\cdots,i_t\in\overline{1,2\m}\}\subset\la
x_n^{\m-k}\theta_1\cdots\theta_{m_1}\ra
\end{equation}
when $k<m_1$. Since
\begin{equation}
x^\alpha\theta_{i_1}\cdots\theta_{i_t}\theta_{2m_1+1}=\frac{(-1)^t}{\alpha_{n+i}+1}(E_{2m+1+1,2m_1+n+i}+E_{2m_1+1+i,2m_1+1})(x_{n+i}x^\alpha\theta_{i_1}\cdots\theta_{i_t}),
\label{b261}
\end{equation}
we obtain ${\cal A}'_k=\la x_n^{m_1-k}\theta_1\cdots\theta_{m_1}\ra$
when $k<\m$. Now we assume $k=m_1$. Since
\begin{eqnarray}
&&(E_{2m_1+1,2\m+1+n}-E_{2m_1+1+2n,2m_1+1})(E_{2m_1+1,2\m+n}-E_{2m_1+2n,2m_1+1})(\theta_1\cdots\theta_{m_1})\nonumber\\
&=&x_{n-1}x_{2n}\theta_1\cdots\theta_{m_1}
\end{eqnarray}
and
\begin{eqnarray}
&&(E_{2m_1+1,2\m+n}-E_{2m_1+2n,2m_1+1})(E_{2m_1+1,2\m+1+n}-E_{2m_1+1+2n,2m_1+1})(\theta_1\cdots\theta_{m_1})\nonumber\\
&=&x_{n}x_{2n-1}\theta_1\cdots\theta_{m_1},
\end{eqnarray}
we have
\begin{equation}
\mbox{Span}\;\{x^\alpha\theta_{i_1}\cdots\theta_{i_t}\in{\cal
A}'_{m_1}\mid 0\leq
t\leq2\m;\;i_1,\cdots,i_t\in\overline{1,2\m}\}\subset\la
x_n^{\m-k}\theta_1\cdots\theta_{m_1}\ra.
\end{equation}
Again by (\ref{b261}), we get ${\cal A}'_{m_1}=\la
\theta_1\cdots\theta_{m_1}\ra$.

Now for any $0\neq f\in{\cal A}'_k$, we can write
$f=f_0+f_1\theta_{2\m+1}$ with $f_0,f_1$ independent of
$\theta_{2m_1+1}$. Applying $osp(2\m,2n)^+$ to $f$, we get some
$0\neq g=g_0+g_1\theta_{2\m+1}\in\la f\ra$ satisfying
$osp(2\m,2n)^+(g_0)=0$ and $osp(2\m,2n)^+(g_1)=0$.

(i) $k<m_1$. According to Theorem 3.2, we can write
\begin{equation}
g=ax_n^{m_1-k}\theta_1\cdots\theta_{m_1}+bx_n^{m_1-k+1}\theta_1\cdots\theta_{m_1}\theta_{2m_1+1}
\end{equation}
where $a,b\in\mathbb{C}$. Thus
$x_n^{m_1-k}\theta_1\cdots\theta_{m_1}\in\la f\ra$ or
$x_n^{m_1-k+1}\theta_1\cdots\theta_{m_1}\theta_{2m_1+1}\in\la f\ra$
because they have different weights. If
$x_n^{m_1-k+1}\theta_1\cdots\theta_{m_1}\theta_{2m_1+1}\in\la f\ra$,
we have
\begin{equation}
x_n^{m_1-k}\theta_1\cdots\theta_{m_1}=\frac{(-1)^{m_1}}{m_1-k+1}(E_{2\m+1,2\m+1+2n}+E_{2\m+1+n,2\m+1})(
x_n^{m_1-k+1}\theta_1\cdots\theta_{m_1}\theta_{2m_1+1}).\label{b262}
\end{equation}
Anyway $x_n^{m_1-k}\theta_1\cdots\theta_{m_1}\in\la f\ra$, which
implies $\la f\ra={\cal A}'_k$.

(ii) $k=m_1$. We can write
\begin{equation}
g=a\theta_1\cdots\theta_{m_1}+b(x_{n-1}x_{2n}-x_nx_{2n-1})\theta_1\cdots\theta_{m_1}+cx_n\theta_1\cdots\theta_{m_1}\theta_{2m_1+1},
\end{equation}
where $a,b,c\in\mathbb{C}$. Thus $\theta_1\cdots\theta_{m_1}\in\la
f\ra$ or
$(x_{n-1}x_{2n}-x_nx_{2n-1})\theta_1\cdots\theta_{m_1}\in\la f\ra$
or $x_n\theta_1\cdots\theta_{m_1}\theta_{2m_1+1}\in\la f\ra$. If
$x_n\theta_1\cdots\theta_{m_1}\theta_{2m_1+1}\in\la f\ra$, we have
$\theta_1\cdots\theta_{m_1}\in\la f\ra$ by (\ref{b262}). If
$(x_{n-1}x_{2n}-x_nx_{2n-1})\theta_1\cdots\theta_{m_1}\in\la f\ra$,
we have
\begin{eqnarray}
&&(E_{2\m+1,2\m+2n}+E_{2\m+n,2\m+1})((x_{n-1}x_{2n}-x_nx_{2n-1})\theta_1\cdots\theta_{m_1})\nonumber\\
&=&-(-1)^{m_1}x_n\theta_1\cdots\theta_{m_1}\theta_{2m_1+1}.
\end{eqnarray}
Thus $\theta_1\cdots\theta_{m_1}\in\la f\ra$.
\end{proof}

\begin{center}{\Large \bf Acknowledgement}\end{center}

I would like to thank Professor Xiaoping Xu for his guidance.

\vspace{1cm}

\bibliographystyle{amsplain}

\end{document}